\documentclass[a4paper, 11pt]{article}
\usepackage{graphicx}
\usepackage{amsmath, amssymb, amsthm}
\usepackage{subfig}
\usepackage{xcolor}
\usepackage{mathrsfs}
\usepackage{multirow}
\usepackage[normalem]{ulem}

\newtheorem{theorem}{Theorem}
\newtheorem{lemma}[theorem]{Lemma}

\newtheorem{definition}{Definition}

\newcommand\bx{\boldsymbol{x}}
\newcommand\bn{\boldsymbol{n}}
\newcommand\bu{\boldsymbol{u}}
\newcommand\bv{\boldsymbol{v}}

\newcommand\bg{\boldsymbol{g}}
\newcommand\bs{\boldsymbol{s}}
\newcommand\bh{\boldsymbol{h}}
\newcommand\bbR{\mathbb{R}}
\newcommand\bbN{\mathbb{N}}
\newcommand\bbS{\mathbb{S}}
\newcommand\czeta{\zeta_{k_1'k_2'k_3',l_1'l_2'l_3'}^{k_1 k_2 k_3,l_1 l_2 l_3}}

\newcommand\dd{\,\mathrm{d}}

\newcommand\mQ{\mathcal{Q}}
\newcommand\mM{\mathcal{M}}
\newcommand\mF{\mathcal{F}}
\newcommand\mI{\mathcal{I}}
\newcommand\bw{\boldsymbol{w}}
\newcommand{\indexk}{k_1k_2k_3}
\newcommand{\indexm}{m_1m_2m_3}
\newcommand{\indexn}{n_1n_2n_3}
\newcommand{\indexl}{l_1l_2l_3}
\newcommand{\newindex}[1]{{{#1}_1{#1}_2{#1}_3}}
\newcommand{\factorialk}{k_1!k_2!k_3!}

\newcommand{\factorialm}{m_1!m_2!m_3!}

\numberwithin{equation}{section}

\graphicspath{{images/}}

{\theoremstyle{remark} }
\newtheorem{corollary}{Corollary}

\title{Approximation of the Boltzmann Collision Operator
  Based on Hermite Spectral Method}

\author{Yanli Wang\thanks{Department of Engineering, Peking University,
  Beijing, China, 100871, email: {\tt wang\_yanli@pku.edu.cn}.}, ~~
  Zhenning Cai\thanks{Department of Mathematics, National
    University of Singapore, Level 4, Block S17, 10 Lower Kent Ridge
    Road, Singapore 119076, email: {\tt matcz@nus.edu.sg}.}}

\begin{document}
\maketitle
\begin{abstract}
Based on the Hermite expansion of the distribution function, we
introduce a Galerkin spectral method for the spatially homogeneous
Boltzmann equation with the realistic inverse-power-law models. A
practical algorithm is proposed to evaluate the coefficients in the
spectral method with high accuracy, and these coefficients are also
used to construct new computationally affordable collision models.
Numerical experiments show that our method captures the low-order
moments very efficiently.

\vspace{3pt}

\end{abstract}

\section{Introduction}
Over a century ago, Boltzmann devised a profound equation describing
the statistical behavior of gas molecules. A number of interesting
theoretical and practical problems emerged due to the birth of this
equation, among which the numerical simulation for this
six-dimensional Boltzmann equation is one significant topic after the
invention of computers. The difficulty comes partly from its high
dimensionality, and partly from its complicated integral operator
modeling the binary collision of gas molecules. People have been using
the Monte Carlo method \cite{Bird} to overcome the difficulty caused
by high dimensionality, but nowadays, a six-dimensional simulation
using a deterministic solver is no longer unaffordable due to the fast
growth of computer flops. Fully six-dimensional computations are
carried out in \cite{Xu2012,Dimarco2015} for the simplified BGK-type
collision terms. However, numerical simulation of the original
Boltzmann equation with the binary collision operator still requires a
large amount of computational resources \cite{Dimarco2018}.

Currently, the deterministic discretization of the binary collision operator
can be categorized into three types: the discrete velocity method
\cite{goldstein1989, Panferov2002}, the Fourier spectral method
\cite{Pareschi1996, Hu2017}, and the Hermite spectral method \cite{Grohs2017,
Gamba2018}. The discrete velocity method is hardly used in the numerical
simulation due to its low order of convergence \cite{Panferov2002}, whereas the
Fourier spectral method is more popular because of its fast convergence rate
and high numerical efficiency. For hard-sphere gases, the computational
cost can be reduced to $O(K N \log N)$ \cite{Mouhot2006, Filbet2006}, with $K$
being the number of discrete points on the unit sphere and $N$ being the total
number of modes in the velocity space. For general gas molecules, the
Fourier spectral method has time complexity $O(N^2)$ \cite{Pareschi1996} or
$O(K N^{4/3} \log N)$ \cite{Hu2017}. Based on these works, some improved
versions of the Fourier spectral methods have been proposed \cite{Gamba2009,
Pareschi2000, Filbet2015, Cai2018}, and some spatially inhomogeneous
applications have been carried out \cite{Wu2014, Dimarco2018}. We would also
like to refer the readers to the review article \cite{Dimarco2014} for a
complete review of the above methods.

Compared with the Fourier spectral method which requires periodization
of the distribution function, the Hermite spectral method looks more
natural since the basis functions with orthogonality in $\mathbb{R}^3$
are employed. In fact, the Hermite spectral method has a longer
history and has been known as the moment method since Grad's work
\cite{Grad}.  Grad proposed in \cite{Grad} a general method to find
the expansion of the binary collision term with Hermite basis
functions. Later, a similar way to expand the binary collision term
using Sonine polynomials (also known as spherical Hermite polynomials)
was proposed in \cite{Kumar}. The techniques used in formulating the
expansion are also introduced in the book \cite{Truesdell}.

Despite these works, the Hermite spectral method is used in the numerical
simulation only until recently \cite{Grohs2017, Gamba2018,
Kitzler2019}. There are two major difficulties in applying this method: one is
the evaluation of the coefficients in the expansion of the collision operator;
the other is the huge computational cost due to its quadratic form. Although
the general procedure to obtain the coefficients is given in \cite{Grad,
Truesdell}, following such a procedure involves expansion of a large number of
huge polynomials, which is quite expensive even for a modern computer algebraic
system; Kumar \cite{Kumar} provided a formula in his expansion using Sonine
polynomials, while the formula involves evaluation of a large number of Talmi
coefficients, which is not tractable either. As for the computational cost, the
computational time of one evaluation of the collision operator is proportional
to the cube of the number of degrees of freedom, while in the Fourier spectral
method, the time complexity for a direct Galerkin discretization is only the
square of the number of modes.

This work is devoted to both of the aforementioned issues. On one
hand, by using a number of properties for relavant polynomials, we
provide explicit formulas for all the coefficients appearing in the
expansion of the collision operator with the Hermite spectral
method. These formulas are immediately applicable in the sense of
coding, and the computational cost is affordable for a moderate number
of degrees of freedom. On the other hand, we combine the modeling
strategy and the numerical technique to form a new way to discretize
the collision term, where only a portion in the truncated series
expansion is treated ``quadratically'', and the remaining part just
decays exponentially as in the BGK model. Thus the computational cost
is greatly reduced and we can still capture the evolution of low-order
moments accurately.

The rest of this paper is organized as follows. In Section
\ref{sec:Bol_Her}, we briefly review the Boltzmann equation and the
Hermite expansion of the distribution function. In Section
\ref{sec:general}, we first give an explicit expression for the series
expansion of the quadratic collision operator, and then construct
approximate collision models based on such an expansion. Some
numerical experiments verifying our method are carried out in Section
\ref{sec:numerical}. Some concluding remarks, as well as some
comparison with similar works, are made in Section
\ref{sec:conclusion}.  Detailed derivation of the expansion is given
in the Appendix.

\section{Boltzmann equation and Hermite expansion of the distribution
  function}
\label{sec:Bol_Her}
This section is devoted to the introduction of existing works needed
by our further derivation. We will first give a brief review of the
Boltzmann equation and the IPL (Inverse-Power-Law) model, and then
introduce the expansion of the distribution function used in the
Hermite spectral method.

\subsection{Boltzmann equation}
\label{sec:Boltzmann}
The Boltzmann equation describes the fluid state using a distribution
function $f(t,\bx,\bv)$, where $t$ is the time, $\bx$ is the spatial
coordinates, and $\bv$ stands for the velocity of gas molecules. The
governing equation of $f$ is 
\begin{equation}
\frac{\partial f}{\partial t} + \nabla_{\bx} \cdot (\bv f)  = \mQ[f],
  \qquad t\in \mathbb{R}^+, \quad \bx \in \mathbb{R}^3,
    \quad \bv \in \mathbb{R}^3,
\end{equation}
where $\mQ[f]$ is the collision operator which has a quadratic form
\begin{equation} \label{eq:quad_col}
  \mQ[f](t,\bx,\bv) = \int_{\mathbb{R}^3}
  \int_{\bn \perp \bg} \int_0^{\pi}
    [f(t,\bx,\bv_1') f(t,\bx,\bv') - f(t,\bx,\bv_1) f(t,\bx,\bv)]
    B(|\bg|,\chi)
  \dd\chi \dd\bn \dd\bv_1,
\end{equation}
where $\bg = \bv - \bv_1$ and $\bn$ is a unit vector. Hence
$\int_{\bn\perp\bg} \cdots \dd\bn$ is a one-dimensional integration
over the unit circle perpendicular to $\bg$. The post-collisional
velocities $\bv'$ and $\bv_1'$ are
\begin{equation} \label{eq:post_vel}
\begin{aligned}
\bv' &= \cos^2(\chi/2) \bv + \sin^2(\chi/2) \bv_1
  - |\bg| \cos(\chi/2) \sin(\chi/2) \bn, \\
\bv_1' &= \cos^2(\chi/2) \bv_1 + \sin^2(\chi/2) \bv
  + |\bg| \cos(\chi/2) \sin(\chi/2) \bn,
\end{aligned}
\end{equation}
and from the conservation of  momentum and energy, it holds that 
\begin{equation}
  \label{eq:mom_eng}
\bv + \bv_1 = \bv' + \bv'_1, \qquad |\bv|^2 + |\bv_1|^2 = |\bv'|^2 +
|\bv'_1|^2. 
\end{equation}
The collision kernel $B(|\bg|, \chi)$ is a non-negative function
determined by the force between gas molecules.

In this paper, we are mainly concerned with the IPL model, for which
the force between two molecules is always repulsive and proportional
to a negative power of their distance.  In this case, the kernel
$B(|\bg|, \chi)$ in \eqref{eq:quad_col} has the form
\begin{equation}
  \label{eq:IPL}
  B(|\bg|, \chi) := |\bg|^{\frac{\eta - 5}{\eta-1}} W_0
  \left| \frac{\mathrm{d} W_0}{\mathrm{d} \chi} \right|, \quad \eta > 3,
\end{equation}
where $-\eta$ is the index in the power of distance. The case $\eta >
5$ corresponds to the ``hard potential'', and the case $3 < \eta < 5$
corresponds to the ``soft potential''. When $\eta = 5$, the collision
kernel $B(|\bg|, \chi)$ is independent of $|\bg|$, and in this model
the gas molecules are called ``Maxwell molecules''. The dimensionless
impact parameter $W_0$ is related to the angle $\chi$ by
\begin{equation} \label{eq:chi}
\chi = \pi - 2 \int_0^{W_1} \left[
  1 - W^2 - \frac{2}{\eta - 1} \left( \frac{W}{W_0} \right)^{\eta-1}
\right]^{-1/2} \dd W,
\end{equation}
and $W_1$ is a positive real number satisfying
\begin{equation} \label{eq:W_1}
1-W_1^2 - \frac{2}{\eta - 1} \left( \frac{W_1}{W_0} \right)^{\eta-1} = 0.
\end{equation}
It can be easily shown that the above equation of $W_1$ admits a
unique positive solution when $\eta > 3$ and $W_0 > 0$.

Apparently, the quadratic collision term is the most complicated part
in the Boltzmann equation. In this paper, we will focus on the
numerical approximation of $\mQ[f]$. For simplicity, we assume that
the gas is homogeneous in space, and thus we can remove the variable
$\bx$ in the distribution function to get the spatially homogeneous
Boltzmann equation
\begin{equation}
  \label{eq:homo}
  \frac{\partial f}{\partial t}   = \mQ[f],
  \qquad t\in \mathbb{R}^+, \quad \bv \in \mathbb{R}^3.
\end{equation}
It is well known that the steady state of this equation takes the form
of the Maxwellian:
\begin{equation}
  \label{eq:general_Maxwellian}
  f(\infty, \bv) = \mathcal{M}_{\rho,\bu,\theta}(\bv)
    := \frac{\rho}{(2\pi\theta)^{3/2}}
      \exp \left( -\frac{|\bv - \bu|^2}{2\theta} \right),
\end{equation}
where the density $\rho$, velocity $\bu$ and temperature $\theta$ can
be obtained by
\begin{equation}
  \label{eq:macro_var}
  \rho = \int_{\bbR^3} f(t,\bv) \dd \bv, \quad \bu =
  \frac{1}{\rho}\int_{\bbR^3}\bv f(t,\bv) \dd \bv, \quad \theta =
  \frac{1}{3\rho}\int_{\bbR^3}|\bv -\bu|^2 f(t,\bv) \dd \bv. 
\end{equation}
These quantities are invariant during the evolution, and therefore
\eqref{eq:macro_var} holds for any $t$. By selecting proper frame of
reference and applying appropriate non-dimensionalization, we can
obtain
\begin{equation}
  \label{eq:u_theta}
  \rho = 1, \quad \bu  = 0, \quad \theta = 1,
\end{equation}
and thus the Maxwellian \eqref{eq:general_Maxwellian} is reduced to
\begin{equation}
  \label{eq:Maxwellian}
  \mathcal{M}(\bv) := \frac{1}{(2\pi)^{3/2}}
    \exp \left( -\frac{|\bv|^2}{2} \right).
\end{equation}
Hereafter, the normalization \eqref{eq:u_theta} will always be
assumed.

In the literature, people have been trying to avoid the complicated
form of the collision operator $\mQ[f]$ by introducing simpler
approximations to it. For example, the BGK collision model
\begin{equation} \label{eq:BGK}
\mQ^{\mathrm{BGK}}[f] =
\frac{1}{\tau} (\mathcal{M} - f)
\end{equation}
was proposed in \cite{BGK}. Here $\tau$ is the mean relaxation time,
which is usually obtained from the first approximation of the
Chapman-Enskog theory \cite{Cowling}. When \eqref{eq:BGK} is used to
approximate the IPL model,
\begin{equation} 
  \label{eq:tau} 
  \tau = \frac{5}{2^{\frac{3\eta-7}{\eta-1}}
    \sqrt{\pi} A_2(\eta)\Gamma(4-2/(\eta-1))},
\end{equation}
where $A_2(\eta) = \int_0^{+\infty} W_0 \sin^2 \chi \dd W_0$. With
$\mQ[f]$ replaced by $\mQ^{\mathrm{BGK}}[f]$ in \eqref{eq:homo}, the
collision process becomes an exponential convergence to the
Maxwellian. Such a simple approximation provides incorrect Prandtl
number $1$. Hence some other models such as the Shakhov model
\cite{Shakhov} and ES-BGK model \cite{Holway} are later proposed to
fix the Prandtl number by changing the Maxwellian in \eqref{eq:BGK} to
a non-equilibrium distribution function. We will call these models
``BGK-type models'' hereafter.

Numerical evaluations on these BGK-type models can be found in
\cite{Gao, chen2015}, where one can find that these approximations are
not accurate enough when the non-equilibrium is strong. Hence the
study on efficient numerical methods for the original Boltzmann
equation with the quadratic collision operator is still necessary.

\subsection{Series expansion of the distribution function}
\label{sec:expansion}
Our numerical discretization will be based on the following series
expansion of the distribution function in the weighted $L^2$ space
$\mF = L^2(\bbR^3; \mM^{-1}\dd\bv)$:
\begin{equation} 
  \label{eq:expansion}
  f(t, \bv) = \sum_{\indexk}f_{k_1k_2k_3}(t)H^{k_1k_2k_3}(\bv)\mM(\bv), 
\end{equation}
where $\mM(\bv)$ is the Maxwellian, and we have used the abbreviation
\begin{equation}
  \label{eq:sumnot}
  \sum_{\indexk} :=  \sum_{k_1=0}^{+\infty} \sum_{k_2=0}^{+\infty}
   \sum_{k_3=0}^{+\infty}.
  \end{equation}
In \eqref{eq:expansion}, $H^{k_1k_2k_3}(\bv)$ are the Hermite
polynomials defined as follows:
\begin{definition}[Hermite polynomials]
  For $k_1, k_2, k_3 \in \mathbb{N}$, the Hermite polynomial $H^{k_1
  k_2 k_3}(\bv)$ is defined as
  \begin{equation}
    \label{eq:basis}
     H^{k_1 k_2 k_3}(\bv) = \frac{(-1)^n}{\mathcal{M}(\bv)}
     \frac{\partial^{k_1+k_2+k_3}}{\partial v_1^{k_1} \partial
       v_2^{k_2} \partial v_3^{k_2}} \mathcal{M}(\bv),
  \end{equation}
  where $\mathcal{M}(\bv)$ is given in \eqref{eq:Maxwellian}.
\end{definition}
The expansion \eqref{sec:expansion} was proposed by Grad in
\cite{Grad}, where such an expansion was used to derive moment
methods. The relation between the coefficients $f_{k_1 k_2 k_3}$ and 
the moments can be seen from the orthogonality of Hermite polynomials
\begin{equation}
  \label{eq:Her_orth}
  \int_{\bbR^3}H^{k_1k_2k_3}(\bv)H^{\indexl}(\bv)\mM(\bv)  \dd \bv
  = \delta_{k_1l_1}\delta_{k_2l_2}\delta_{k_3l_3}k_1!k_2!k_3!.
\end{equation}
For example, by the above orthogonality, we can insert the expansion
\eqref{eq:expansion} into the definition of $\rho$ in
\eqref{eq:macro_var} to get $f_{000} = \rho$. In our case, the
normalization \eqref{eq:u_theta} gives us $f_{000} = 1$. Similarly,
it can be deduced from the other two equations in \eqref{eq:macro_var}
and \eqref{eq:u_theta} that
\begin{equation}
  f_{100} = f_{010} = f_{001} = 0, \qquad f_{200} + f_{020} + f_{002} = 0.
\end{equation}
Other interesting moments include the heat flux $q_i$ and the stress
tensor $\sigma_{ij}$, which are defined as
\begin{gather*}
  q_i = \frac{1}{2}\int_{\bbR^3}|\bv -\bu|^2(v_i-u_i)f\dd \bv 
  = \frac{1}{2}\int_{\bbR^3}|\bv|^2v_if\dd \bv, \qquad i = 1, 2, 3, \\
\sigma_{ij} = \int_{\bbR^3}
  \left( (v_i-u_i)(v_j -u_j) - \frac{1}{3} \delta_{ij} |\bv-\bu|^2 \right) f
\dd \bv =\int_{\bbR^3}
  \left( v_iv_j - \frac{1}{3} \delta_{ij} |\bv|^2 \right) f
\dd \bv, \quad i,j = 1,2,3. 
\end{gather*}
They are related to the coefficients by
\begin{displaymath}
q_1 =  3f_{300} + f_{120} + f_{102}, \qquad
q_2 = 3f_{030} + f_{210} + f_{012}, \qquad
q_3 = 3f_{003} + f_{201} + f_{021},
\end{displaymath}
and
\begin{gather*}
\sigma_{11} = 2f_{200}, \qquad
\sigma_{12} = f_{110}, \qquad
\sigma_{13} = f_{101}, \\
\sigma_{22} = 2f_{020}, \qquad
\sigma_{23} = f_{011}, \qquad
\sigma_{33} = 2 f_{002}.
\end{gather*}


\section{Approximation of the collision term}
\label{sec:general}
To get the evolution of the coefficients $f_{k_1 k_2 k_3}$ in the
expansion \eqref{eq:expansion}, we need to expand the collision term
using the same basis functions. The expansions of the BGK-type
collision operators are usually quite straightforward. For instance,
the series expansion of the BGK collision term \eqref{eq:BGK} is given
in \cite{Cai} as
\begin{equation}
  \label{eq:BGK_exp}
  \mathcal{Q}^{\mathrm{BGK}}[f] = 
  \sum_{\indexk}Q_{\indexk}^{\mathrm{BGK}} H^{k_1k_2k_3}(\bv)\mM(\bv),
\end{equation}
where 
\begin{equation*}
  Q_{\indexk}^{\mathrm{BGK}} = \left\{
      \begin{array}{cl}
        0, &  k_1 = k_2 = k_3 = 0, \\
        -\frac{1}{\tau}f_{\indexk}, & \text{otherwise}.
      \end{array}
    \right.
\end{equation*}
The expansions for the ES-BGK and Shakhov operators can be found in
\cite{Cai2013, Cai2012}. In this section, we will first discuss the
series expansion of the quadratic collision term $\mQ[f]$ defined in
\eqref{eq:quad_col}, and then mimic the BGK-type collision operators
to construct collision models with better accuracy.

\subsection{Series expansions of general collision terms}
Suppose the binary collision term $\mQ[f]$ can be expanded as
\begin{equation}
  \label{eq:S_expan}
  \mQ[f](\bv) = \sum_{\indexk}
  Q_{k_1k_2k_3}H^{k_1k_2k_3}(\bv)\mM(\bv).
\end{equation}
By the orthogonality of Hermite polynomials, we get
\begin{equation}
\label{eq:S_k1k2k3}
Q_{\indexk} 
= \frac{1}{k_1!k_2!k_3!}\int H^{k_1k_2k_3}(\bv)\mQ[f](\bv) \dd \bv  
=\sum\limits_{\newindex{i}}\sum\limits_{\newindex{j}} 
  A_{k_1k_2k_3}^{\newindex{i}, \newindex{j}}
  f_{\newindex{i}}f_{\newindex{j}},
\end{equation}
where the second equality can be obtained by inserting
\eqref{eq:expansion} into \eqref{eq:quad_col}, and
\begin{equation}
  \label{eq:coeA_detail}
  \begin{aligned}
    A_{k_1k_2k_3}^{\newindex{i}, \newindex{j}} =&
    \frac{1}{(2\pi)^3\factorialk} \int_{\mathbb{R}^3}
    \int_{\mathbb{R}^3} \int_{\bn \perp \bg} \int_0^{\pi}
    B(|\bg|,\chi)
    \left[H^{\newindex{i}}(\bv')H^{\newindex{j}}(\bv'_1)\right. \\
    & \left. -H^{\newindex{i}}(\bv)H^{\newindex{j}}(\bv_1) \right]
    H^{k_1 k_2 k_3}(\bv) \exp \left( -\frac{|\bv|^2 + |\bv_1|^2}{2}
    \right) \,\mathrm{d}\chi \,\mathrm{d}\bn \,\mathrm{d}\bv_1
    \,\mathrm{d}\bv.
\end{aligned}
\end{equation}
It can be seen from \eqref{eq:coeA_detail} that the evaluation of
every coefficient requires integration of an eight-dimensional
function. In principle, this can be done by numerical quadrature;
however, the computational cost for obtaining all these coefficients
would be huge.  Actually, in \cite{Grad, Truesdell}, a strategy to
simplify the above integral has been introduced, and for small
indices, the values are given in the literature. However, when the
indices are large, no explicit formulae are provided in \cite{Grad,
Truesdell}, and the procedure therein is not easy to follow. Inspired
by these works, we give in this paper explicit equations of the
coefficients $A_{k_1k_2k_3}^{\newindex{i}, \newindex{j}}$ for any
collision kernel, except for an integral with respect the two
parameters in the kernel function $B(\cdot,\cdot)$. The main results
are summarized in the following two theorems:

\begin{theorem} 
  \label{thm:coeA}
  The expansion coefficients of the collision operator $\mQ[f](\bv)$
  defined in \eqref{eq:S_k1k2k3} have the form below:
  \begin{equation}
    \label{eq:coeA}
    \begin{aligned}
      A_{k_1k_2k_3}^{\newindex{i}, \newindex{j}}
      =&  \sum\limits_{i'_1=0}^{\min(i_1+j_1,k_1)}
      \sum\limits_{i'_2=0}^{\min(i_2+j_2,k_3)}
      \sum\limits_{i'_3=0}^{\min(i_3+j_3,k_3)}
      \frac{ 2^{-k/2}}{2^3\pi^{3/2}}\frac{1}{l'_1!l'_2!l'_3!}
      a_{i'_1j'_1}^{i_1j_1}a_{i'_2j'_2}^{i_2j_2}
      a_{i'_3j'_3}^{i_3j_3}\gamma_{j'_1j'_2j'_3}^{l'_1l'_2l'_3},
    \end{aligned}
  \end{equation}
  where
  \begin{equation}
    \label{eq:ijl_A}
    j'_s = i_s+j_s-i'_s,\quad l'_s=k_s -i'_s, \qquad s = 1, 2, 3.    
  \end{equation}
  The coefficients $a_{i'j'}^{ij}$ and
  $\gamma_{j'_1j'_2j'_3}^{l'_1l'_2l'_3}$ are defined by
  \begin{equation}
    \label{eq:coea}
    a_{i'j'}^{ij} = 2^{-(i'+j')/2} i! j!
    \sum_{s=\max(0,i'-j)}^{\min(i',i)}
    \frac{(-1)^{j'-i+s}}{s!(i-s)!(i'-s)!(j'-i+s)!},
  \end{equation}
  and
  \begin{equation}
    \label{eq:coe_gamma}
    \gamma_{\newindex{j}}^{\indexl} :=
    \int_{\mathbb{R}^3} \int_{\bn \perp \bg} \int_0^{\pi} \left[
      H^{\newindex{j}} \left( \frac{\bg'}{\sqrt{2}} \right) -
      H^{\newindex{j}} \left( \frac{\bg}{\sqrt{2}} \right)
    \right]
    H^{\indexl} \left( \frac{\bg}{\sqrt{2}} \right)
    B(|\bg|, \chi) \exp \left( -\frac{|\bg|^2}{4} \right)
    \,\mathrm{d}\chi \,\mathrm{d}\bn \,\mathrm{d}\bg,
  \end{equation}
  where $\bg' = \bg\cos\chi - |\bg| \bn \sin\chi$ is the
  post-collisional relative velocity, and $B(|\bg|, \chi)$ is the
  collision kernel in \eqref{eq:quad_col}.
\end{theorem}
\begin{theorem}
  \label{thm:gamma}
  For any $k_1, k_2, k_3, l_1, l_2, l_3 \in \mathbb{N}$, let $k = k_1
  + k_2 + k_3$ and $l = l_1 + l_2 + l_3$. Then the coefficients
  $\gamma_{{\newindex{j}}}^{\indexl}$ defined in \eqref{eq:coe_gamma}
  satisfies
  \begin{equation}
    \label{eq:gamma}
    \gamma_{\indexk}^{l_1 l_2 l_3} = 
    \sum_{m_1=0}^{\lfloor k_1/2 \rfloor}
    \sum_{m_2=0}^{\lfloor k_2/2 \rfloor}
    \sum_{m_3=0}^{\lfloor k_3/2 \rfloor}
    \sum_{n_1=0}^{\lfloor l_1/2 \rfloor}
    \sum_{n_2=0}^{\lfloor l_2/2 \rfloor}
    \sum_{n_3=0}^{\lfloor l_3/2 \rfloor}(2k-4m+1)
    C_{m_1 m_2 m_3}^{\indexk} C_{\indexn}^{l_1 l_2 l_3}
    S_{k_1-2m_1,k_2-2m_2,k_3-2m_3}^{l_1-2n_1,l_2-2n_2,l_3-2n_3}
    K_{mn}^{kl},
  \end{equation}
  where $m = m_1 + m_2 + m_3$, $n = n_1 + n_2 + n_3$, and
  \begin{equation}
    \label{eq:coeC}
    C_{\indexm}^{\indexk} =
    \frac{(-1)^m4\pi m!} {(2(k-m)+1)!!}\frac{\factorialk}{\factorialm}.
  \end{equation}
  In \eqref{eq:gamma}, $S_{\indexk}^{\indexl}$ is the coefficient of
  $v_1^{k_1} v_2^{k_2} v_3^{k_3} w_1^{l_1} w_2^{l_2} w_3^{l_3}$ in the
  polynomial
  \begin{equation}
    \label{eq:Sk}
    S_k(\bv, \bw) := (|\bv| |\bw|)^k P_k \left(
      \frac{\bv}{|\bv|} \cdot \frac{\bw}{|\bw|}
    \right), 
  \end{equation}
  and 
  \begin{equation}
    \label{eq:Klmn}
    \begin{split}
      K_{mn}^{kl} & = \int_0^{+\infty} \int_0^{\pi}
      L_m^{(k-2m+1/2)}\left( \frac{g^2}{4} \right)
      L_n^{(l-2n+1/2)}\left( \frac{g^2}{4} \right) \\
      & \qquad \times \left( \frac{g}{\sqrt{2}} \right)^{k+l+2-2(m+n)}
      B(g,\chi) \Big[P_{k-2m}(\cos \chi) - 1\Big]
      \exp \left( -\frac{g^2}{4} \right) \,\mathrm{d}\chi \,\mathrm{d}g.
    \end{split}
  \end{equation}
  Here $L_n^{(\alpha)}(x)$ are the Laguerre polynomials and $P_k(x)$ are
  the Legendre polynomials, which are defined below.
\end{theorem}

\begin{definition}[Legendre functions]
  For $\ell \in \mathbb{N}$, the Legendre polynomial $P_{\ell}(x)$ is
  defined as
\begin{displaymath}
  P_{\ell}(x) = \frac{1}{2^{\ell} \ell!}
  \frac{\mathrm{d}^{\ell}}{\mathrm{d} x^{\ell}}
  \left[ (x^2 - 1)^{\ell} \right].
\end{displaymath}
\end{definition}

\begin{definition}[Laguerre polynomials]
\label{def:Laguerre}
For $\alpha > -1$, let $w^{\alpha}(x) = x^{n+\alpha} \exp(-x)$. For $n
\in \mathbb{N}$, define the Laguerre polynomial as
\begin{displaymath}
L_n^{(\alpha)}(x) = \frac{x^n}{n! w^{\alpha}(x)}
  \frac{\mathrm{d}^n}{\mathrm{d}x^n} w^{\alpha}(x).
\end{displaymath}
\end{definition}

Through these two theorems, the eight-dimensional integration in
\eqref{eq:coeA_detail} has been reduced into a series of summations
and a two-dimensional integration. Among all the coefficients
introduced in these theorems, $a_{i'j'}^{ij}$ and
$C_{\indexm}^{\indexk}$ can be computed directly. As for
$S_{\indexk}^{\indexl}$, we need to expand polynomial $S_k(\bv, \bw)$,
which can be done recursively using the following recursion formula:
\begin{equation}
  \label{eq:recu_Sk}
  \begin{gathered}
    S_0(\bv, \bw) = 1, \qquad S_1(\bv, \bw) = \bv\cdot \bw,\\
    S_{k+1}(\bv, \bw) = \frac{2k+1}{k+1}(\bv\cdot \bw)S_k(\bv, \bw) -
    \frac{k}{k+1}(|\bv||\bw|)^2S_{k-1}(\bv, \bw).
  \end{gathered}
\end{equation} 
This recursion formula can be derived from the recursion relation of
Legendre polynomials, and it also shows that for every monomial in the
expansion of $S_k(\bv, \bw)$, the degree of $\bv$ equals the degree of
$\bw$. Therefore $S_{\indexk}^{\indexl}$ is nonzero only when $k_1 +
k_2 + k_3 = l_1 + l_2 + l_3$. This means in \eqref{eq:gamma}, the
summand is nonzero only when 
\begin{equation}
  \label{eq:kl_gamma}
  k_1 + k_2 + k_3 - 2(m_1 + m_2 + m_3) =
    l_1 + l_2 + l_3 - 2(n_1 + n_2 + n_3).
\end{equation}
Consequently, when evaluating $K_{mn}^{kl}$ defined in
\eqref{eq:Klmn}, we only need to take into account the case $k-2m =
l-2n$. Generally, $K_{mn}^{kl}$ can be computed by numerical
quadrature; for the IPL model, the integral with respect to $g$ can
be written explicitly, which will be elaborated in the following
section.

\subsection{Series expansion of collision operators for IPL models}
\label{sec:IPL}
The formulae given in the previous section are almost ready to be
coded, except that specific collision models are needed to calculate
the integral $K_{mn}^{kl}$ defined in \eqref{eq:Klmn}. This section is
devoted to further simplifying this integral for IPL models, which
completes the algorithm for computing the coefficients
$A_{\indexk}^{\newindex{i}, \newindex{j}}$.

For the IPL model \eqref{eq:IPL}, we first consider the integral with
respect to $\chi$ in \eqref{eq:Klmn}. To this aim, we extract all the
terms related to $\chi$ from \eqref{eq:Klmn}, and define
$\tilde{B}_{k}^{\eta}(\cdot)$ as
\begin{equation*}
  \tilde{B}_{k}^{\eta}(g) := \int_0^{\pi}
  B(g,\chi) \Big[P_{k}(\cos \chi) - 1\Big]
  \dd \chi = g^{\frac{\eta - 5}{\eta-1}}\int_0^{\pi} W_0
  \left| \frac{\mathrm{d} W_0}{\mathrm{d} \chi} \right|
  \Big[P_{k}(\cos \chi) - 1\Big]
  \dd \chi, \qquad \eta>3, \quad g > 0.
\end{equation*} 
To evaluate the above integral, we follow the method introduced in
\cite{Cai2015} and apply the change of variable
\begin{displaymath}
  \chi = \pi - 2\int_0^1 [1-x^2(1-y)-x^{\eta-1} y]^{-1/2} \sqrt{1-y} \dd x,
\end{displaymath}
to get
\begin{equation}
  \label{eq:tB_exact}
  \tilde{B}_{k}^{\eta}(g) = 2^{-\frac{\eta-3}{\eta-1}}
    g^{\frac{\eta - 5}{\eta-1}}\int_0^1 [P_k(\cos \chi) - 1]
  [2(1-y) + (\eta-1)y] [(\eta-1)y]^{-\frac{\eta+1}{\eta-1}} \dd y,
\end{equation}
Below we write the above equation as
\begin{equation}
  \label{eq:tB_result}
  \tilde{B}_{k}^{\eta}(g) = 2^{-\frac{\eta-3}{\eta-1}}
    g^{\frac{\eta - 5}{\eta-1}} \mI(k,\eta),
\end{equation}
where $\mI(k, \eta)$ denotes the integral in \eqref{eq:tB_exact}. In
general, we need to evaluate $\mI(k, \eta)$ by numerical quadrature.
In our implementation, the adaptive integrator introduced in
\cite[Section 3.3.7]{Piessens} is used to compute this integral.

Now we consider the integral with respect to $g$. Using the result
\eqref{eq:tB_result}, we can rewrite \eqref{eq:Klmn} as
\begin{equation}
  \label{eq:Klmn_short}
  K_{mn}^{kl} = 2^{c(\eta)}\mI(k-2m,\eta)\int_{0}^{+\infty} L_m^{(k-2m+1/2)}(s)
  L_n^{(k-2m+1/2)}(s) s^{c(\eta)}\exp(-s) \dd s,
\end{equation}
where $c(\eta) = \frac{\eta-3}{\eta-1} + k - 2m$, and we have applied
the change of variable $s = g^2/4$, and taken into account the
relation $k-2m=l-2n$. In general, we can adopt the formula
\begin{equation} \label{eq:no_hg}
  \begin{split}
    \int_0^{+\infty} L_{m}^{(\alpha)} (s) L_{n}^{(\alpha)} (s) s^{\mu} \exp(-s)
    \dd s = (-1)^{m+n}\Gamma(\mu+1)\sum_{i=0}^{\min(m,n)} \binom{\mu
      - \alpha}{m - i} \binom{\mu - \alpha}{n - i} \binom{i + \mu}{i}
  \end{split}
\end{equation}
introduced in \cite[eq. (10)]{Srivastava} to calculate
\eqref{eq:Klmn_short}. Specially, when $\eta = 5$, which corresponds
to the model of Maxwell molecules, we can use the orthogonality of
Laguerre polynomials to get
\begin{equation} \label{eq:no_hg_delta}
  K_{mn}^{kl} = 2^{k-2m+1/2}\mI(k-2m,5)
    \binom{k-m+1/2}{m} \Gamma(k-2m+3/2) \delta_{mn},
\end{equation}
and thus the computational cost can be further reduced. In fact, Grad
has already pointed out in \cite{Grad} that for Maxwell molecules,
$A_{\indexk}^{\newindex{i}, \newindex{j}}$ is nonzero only when
\begin{equation} \label{eq:i+j=k}
i_1 + i_2 + i_3 + j_1 + j_2 + j_3 = k_1 + k_2 + k_3.
\end{equation}
This can also be seen from our calculation: from
\eqref{eq:no_hg_delta}, we can find that only when $k_1 + k_2 + k_3 =
l_1 + l_2 + l_3$, the coefficient $\gamma_{\indexk}^{\indexl}$ given
in \eqref{eq:gamma} is nonzero; therefore in \eqref{eq:coeA}, if the
summand is nonzero, the sum of $j'_1$, $j'_2$ and $j'_3$ must equal
the sum of $l'_1$, $l'_2$ and $l'_3$, which is equivalent to
\eqref{eq:i+j=k} due to \eqref{eq:ijl_A}.

The above analysis shows that for the IPL model, we only need to apply
the numerical quadrature to the one-dimensional integrals
$\mI(k,\eta)$, which makes it easier to obtain the coefficients
$A_{\indexk}^{\newindex{i},\newindex{j}}$ with high accuracy.

\subsection{Approximation of the collision term}
\label{sec:approx}
Until now, we already have a complete algorithm to calculate the
coefficients $A_{\indexk}^{\newindex{i},\newindex{j}}$. These
coefficients can be used either to discretize the collision term or to
construct new collision models. We will discuss both topics in this
section.

\subsubsection{Discretization of the homogeneous Boltzmann equation}
\label{sec:algo} 
Based on the expansion of the distribution function
\eqref{eq:expansion}, the most natural discretization of the
homogeneous Boltzman equation is to use the Galerkin spectral method.
From this point of view, for any positive integer $M$, we define the
space of the numerical solution
\begin{equation}
  \label{eq:PM}
  \mF_M =  \mathrm{span} \{
  H^{\indexk}(\bv)\mM(\bv) \mid (k_1, k_2,k_3) \in I_M\} \subset \mF =
  L^2(\mathbb{R}^3; \mM^{-1}\dd\bv),
\end{equation}
where $I_M$ is the index set
\begin{displaymath}
I_M = \{(k_1, k_2, k_3) \mid 0\leqslant k_1 + k_2 + k_3 \leqslant M,
  \: k_i \in \bbN, \: i = 1,2,3\}.
\end{displaymath}
Then the semi-discrete distribution function $f_M(t,\cdot) \in \mF_M$
satisfies
\begin{equation} \label{eq:var}
\int_{\mathbb{R}^3} \frac{\partial f_M}{\partial t} \varphi
  \mathcal{M}^{-1} \,\mathrm{d}\bv =
\int_{\mathbb{R}^3} \mQ(f_M, f_M) \varphi \mathcal{M}^{-1} \dd \bv,
  \qquad \forall \varphi \in \mF_M.
\end{equation}
Suppose
\begin{equation} 
  \label{eq:f_h}
  f_M(t,\bv) = \sum\limits_{(k_1,k_2,k_3)\in I_M}
    f_{k_1k_2k_3}(t)H^{k_1k_2k_3}(\bv)\mM(\bv) \in \mF_M.
\end{equation}
The equations \eqref{eq:S_expan} and \eqref{eq:S_k1k2k3} show that the
variational form \eqref{eq:var} is equivalent to the following ODE
system:
\begin{equation} \label{eq:ODE}
\frac{\mathrm{d}f_{k_1k_2k_3}}{\mathrm{d}t} =
  \sum\limits_{(i_1,i_2,i_3)\in I_M} \sum\limits_{(j_1,j_2,j_3)\in I_M} 
    A_{k_1k_2k_3}^{\newindex{i},\newindex{j}}f_{\newindex{i}}f_{\newindex{j}},
\qquad (k_1,k_2,k_3) \in I_M.
\end{equation}
It is easy to see that the time complexity for the computation of all
the right-hand sides is proportional to the number of nonzero
coefficients. For most collision operators, the coefficients
$A_{k_1k_2k_3}^{\newindex{i}, \newindex{j}}$ form a full tensor, since
there is no evidence showing that
$A_{k_1k_2k_3}^{\newindex{i}, \newindex{j}}$ can be zero, except a few
coefficients related to the conservation laws. Therefore, the
computational cost for the right-hand side of \eqref{eq:ODE} is
$O(N_M^3) = O(M^9)$, where $N_M$ is the number of elements in $I_M$:
\begin{equation}
  \label{eq:coeN}
  N_M = \frac{(M+1)(M+2)(M+3)}{6} \sim O(M^3).
\end{equation}
However, when considering Maxwell molecules, due to the
constraint \eqref{eq:i+j=k}, the computational cost can be reduced to
$O(M^8)$.

To fully formulate the ODE system \eqref{eq:ODE}, we need the
coefficients $A_{k_1k_2k_3}^{\newindex{i},\newindex{j}}$ for all
$(i_1,i_2,i_3), (j_1,j_2,j_3), (k_1,k_2,k_3) \in I_M$. When the
collision kernel is chosen and $M$ is fixed, we only need to compute
these coefficients once, and then they can be used repeatedly. For a
given $M$, the algorithm for computing these coefficients is
summarized in Table \ref{tab:index1}. The general procedure is to
sequentially compute the coefficients in the first column, with
indices described in the third column, and the equations to follow are
given in the second column. For IPL models, we can use
\eqref{eq:Klmn_short} and \eqref{eq:no_hg} instead to obtain the
values of $K_{mn}^{kl}$.  In the third column of Table
\ref{tab:index1}, it is worth mentioning that some indices are in the
index set $I_{2M}$ instead of $I_M$, as is due to the equation
\eqref{eq:ijl_A}, which shows that
\begin{displaymath}
(j'_1, j'_2, j'_3) \in I_{2M}, \quad
  \text{if } (i_1, i_2, i_3) \in I_M \text{ and } (j_1, j_2, j_3) \in I_M.
\end{displaymath}
Therefore the corresponding indices for $\gamma$ and $C$ must lie in
$I_{2M}$. Similar arguments hold for the coefficients $K$.

The last column in Table \ref{tab:index1} shows an estimation of the
computational cost for each coefficient, from which one can see that
the total cost for getting $A_{k_1k_2k_3}^{\newindex{i},\newindex{j}}$
is $O(M^{12})$. Now we compare this with the numerical cost by
applying numerical integration directly to \eqref{eq:coeA_detail}.
We assume the number of quadrature points on $\mathbb{R}^3$ is
$O(M_v^3)$, and the number of quadrature points on the unit sphere
(domain for $\bn$ and $\chi$) is $O(M_s^2)$. Thus using numerical
integration to evaluate all the coefficients
$A_{k_1k_2k_3}^{\newindex{i},\newindex{j}}$ has time complexity $O(M^9
M_v^6 M_s^2)$. In most cases, we will choose $M_v > M$ to get accurate
results. Hence our method listed in Table \ref{tab:index1} is
significantly faster.

\begin{table}[hpt!]
  \centering
  \def\arraystretch{1.5}
  {\footnotesize
  \begin{tabular}{cccc}
    \hline
    Coefficients & Formula & Constraints for the indices & Computational cost\\
    \hline 
    $C_{\indexm}^{\indexk}$ & \eqref{eq:coeC} 
                           &  $(k_1,k_2, k_3) \in I_{2M},$ \quad $(m_1,m_2, m_3) \in I_{M}$  & $O(M^6)$\\
    $S_{\indexl}^{\indexk}$ & \eqref{eq:recu_Sk} 
                           & $(k_1, k_2,k_3) \in I_{M}$, \quad
			     $(l_1, l_2,l_3) \in I_{M}$, \quad $k_1 + k_2 + k_3 = l_1 + l_2 + l_3$ & $O(M^5)$\\
    $K_{mn}^{kl}$ &  \eqref{eq:Klmn}  
                           & $k \leqslant 2M$, \quad
                               $l \leqslant M$, \quad
                               $m\leqslant \lfloor k/2 \rfloor$, \quad
                               $n\leqslant \lfloor l/2 \rfloor$, \quad
                               $k -2m = l - 2n$ & $O(M^4)$\\
    $\gamma_{\indexk}^{\indexl}$ & \eqref{eq:gamma} 
                           & $(l_1,l_2,l_3)\in I_{M}$, \quad $(k_1, k_2, k_3) \in I_{2M}$ &$O(M^{11})$ \\
    $a_{i'j'}^{ij}$ & \eqref{eq:coea} 
                           & $i \leqslant  M$, \quad $j \leqslant  M$,\quad
                             $i' \leqslant  2M$, \quad $j' \leqslant 2M$, \quad
                             $i + j = i' + j'$ & $O(M^4)$\\
    $A_{\indexk}^{\newindex{i},\newindex{j}}$ & \eqref{eq:coeA} 
                           & $(k_1,k_2,k_3)\in I_M$,\quad $(i_1,i_2, i_3)  \in I_M$, \quad
                             $(j_1, j_2, j_3) \in I_M$  & $O(M^{12})$\\
    \hline
  \end{tabular}}
\caption{A summary for computation of all the coefficients.}
\label{tab:index1}
\end{table}

\subsubsection{Approximation of the collision operator}
\label{sec:convergence}
In the previous section, a complete numerical method has been given to
solve the spatially homogeneous Boltzmann equation. However, due to
the rapid growth of the number of coefficients as $M$ increases, the
storage requirement of this algorithm is quite strong. Table
\ref{tab:memory} shows the memory required to store the coefficients
$A_{\indexk}^{\newindex{i},\newindex{j}}$, where we assume that the
coefficients are represented in the double-precision floating-point
format, whose typical size is $8$ bytes per number. It can be seen
that the case $M = 20$ has already exceeded the memory caps of most
current desktops. Although the data given in Table \ref{tab:memory}
can be reduced by taking the symmetry of the coefficients into
consideration, it can still easily hit our memory limit by increasing
$M$ slightly. Even if the memory cost is acceptable for large $M$, the
computational cost $O(M^9)$ becomes an issue especially when solving
the spatially inhomogeneous problems.
\begin{table}[!ht]
\centering
\begin{tabular}{cccc}
\hline
$M$ & Memory (Gigabytes) & $M$ & Memory (Gigabytes) \\
\hline
5 & $1.308 \times 10^{-3}$ & 25 & $2.620 \times 10^2$ \\
10 & $0.1743$ & 30 & $1.210 \times 10^3$ \\
15 & $4.048$ & 35 & $4.473 \times 10^3$ \\
20 & $41.38$ & 40 & $1.400 \times 10^4$ \\
\hline
\end{tabular}
\caption{Memory required to store $A_{\indexk}^{\newindex{i},\newindex{j}}$.}
\label{tab:memory}
\end{table}

To overcome this difficulty, we will only compute and store the
coefficients $A_{\indexk}^{\newindex{i},\newindex{j}}$ for a small
number $M$ such that the computational cost for solving
\eqref{eq:ODE} is acceptable. When $(k_1, k_2, k_3) \not\in I_M$, we
apply the idea of the BGK-type models and let these coefficients decay
to zero exponentially with a constant rate:
\begin{equation} \label{eq:decay}
\frac{\mathrm{d}f_{\indexk}}{\mathrm{d}t}  = -\nu_M f_{\indexk},
  \qquad (k_1, k_2, k_3) \not\in I_M,
\end{equation}
where $\nu_M$ is a constant independent of $k_1$, $k_2$ and $k_3$.
Combining \eqref{eq:ODE} and \eqref{eq:decay}, we actually get a new
collision operator
\begin{equation}
  \label{eq:approQ}
  \mQ_M[f] = P_M \mQ[P_Mf] - \nu_M(I - P_M)f, \qquad \forall f \in \mF,
\end{equation}
where $P_M$ is the orthogonal projection from $\mF$ onto $\mF_M$.
Such a idea is to mimic the derivation of Shakhov model
\cite{Shakhov}, which models the collision by
\begin{equation}
\label{eq:Shakhov}
\mQ^{\mathrm{S}}[f] := P_{G13} \mathcal{L}[P_{G13} f]
  - \nu(I - P_{G13})f, \qquad \forall f \in \mF,
\end{equation}
where $\mathcal{L}$ is the linearized collision operator defined by
\begin{displaymath}
\mathcal{L}[f] := \lim_{\epsilon \rightarrow 0}
  \frac{Q[\mM + \epsilon(f - \mM)]}{\epsilon}
\end{displaymath}
and $P_{G13}$ is the projection operator onto the 13-dimensional
subspace
\begin{displaymath}
\mF_{G13} = \left\{
  p(\bv) \mM(\bv) \,\Bigg|\, p(\bv) = \alpha +
    \sum_{j=1}^3 \beta_j v_j + \sum_{i,j=1}^3 \gamma_{ij} v_i v_j +
    \sum_{j=1}^3 \zeta_j |\bv|^2 v_j
\right\}
\end{displaymath}
which includes Grad's 13 moments \cite{Grad}. Comparing
\eqref{eq:approQ} and \eqref{eq:Shakhov}, one finds that in our model,
we have replaced the linearized collision operator $\mathcal{L}$ by
the more accurate quadratic collision operator $\mathcal{Q}$, and the
subspace $\mF_{G13}$ is replaced by the larger space $\mF_M$ once $M
\geqslant 3$. Thus the proposed model is expected to provide better
accuracy than the Shakhov model.

The difference between the proposed model and the original quadratic
model is to be further studied in the future work. In general, we
suppose
\begin{enumerate}
\item The projection operator $P_M$ has spectral accuracy;
\item $\mQ[P_M f]$ approximates $\mQ[f]$ with spectral accuracy.
\end{enumerate}
Then
\begin{displaymath}
\|\mQ_M[f] - \mQ[f]\| \leqslant \|P_M \mQ[P_M f] - \mQ[P_M f]\|
  + \|\mQ[P_M f] - \mQ[f]\| + |\nu_M| \|f - P_M f\|,
\end{displaymath}
from which one can see that $\mQ^S[f]$ approximates $\mQ[f]$ with
spectral accuracy. Applying spectral method to this collision operator
is quite straightforward. One just needs to choose an appropriate $M$
(modelling parameter) and an appropriate index set for $k_1, k_2$ and
$k_3$ (discretization parameter), and then solve the ODE system
combined by \eqref{eq:ODE} and \eqref{eq:decay} for $k_1, k_2, k_3$ in
the index set. Thus, it remains only to select the constant $\nu_M$.

In \cite{Cai2015}, the authors used a similar idea to approximate the
linearized collision operator, where the evolution of the coefficients
for high-degree basis functions is also approximated by an exponential
decay. Here we choose the decay rate in the same way as in
\cite{Cai2015}: considering the discrete linearized collision operator
$\mathcal{L}_M: \mF_M \rightarrow \mF_M$ defined as
\begin{equation}
  \mathcal{L}_M[f] =\sum\limits_{(k_1,k_2,k_3) \in I_M}
  \sum\limits_{(j_1, j_2, j_3) \in I_M}
  (A_{\indexk}^{000,\newindex{j}} +  A_{\indexk}^{\newindex{j},000})
  f_{\newindex{j}} H^{\indexk}(\bv) \mM(\bv),
\end{equation}
we let $\nu_M$ be the spectral radius of this operator. The idea of
such a choice includes the following:
\begin{enumerate}
\item As ``less important coefficients'' ($(k_1, k_2, k_3) \not\in
  I_M$), the decay rate should be faster than all the ``important
  coefficients'' ($(k_1, k_2, k_3) \in I_M$). Therefore we choose
  $\nu_M \geqslant \rho(\mathcal{L}_M)$, where $\rho(\mathcal{L}_M)$
  is the spectral radius of $\mathcal{L}_M$, indicating the fastest
  decay rate for the important coefficients.
\item We do not want to introduce any gap between the spectrum of the 
  two parts, causing a sharp transition in the frequency space.
  Therefore we choose $\nu_M = \rho(\mathcal{L}_M)$.
\end{enumerate}
Additionally, it has also been shown in \cite{Cai2015} that such a
choice of $\nu_M$ agrees with the choice of $\nu$ in the Shakhov model
\eqref{eq:Shakhov}. By taking the same $\nu_M$ in $\mQ_M[f]$, the
linearization of $\mQ_M[f]$ about the Maxwellian $\mathcal{M}$
coincides with the approximation of the linearized collision operator
proposed in \cite{Cai2015}.

The collision operator $\mathcal{Q}_M$ deals with a high-frequency
modes with a very simple method: they are damped to zero at a uniform
decay rate. However, in the solution of the Boltzmann equation, it is
often observed that higher-frequency modes decay faster (see Section
\ref{sec:BKW} for an example). This can be achieved by a more careful
modelling for the higher-frequency modes.  Although not yet
implemented, we would like to discuss some possibilities to make
improvements. The first possibility is to replace the simple uniform
decay by the linearized collision operator:
\begin{displaymath}
\mQ^*_M[f] = P_M \mQ[P_M f] + \mathcal{L}[(I-P_M)f].
\end{displaymath}
Since the computation of the linearized collision operator is much
cheaper than that of the quadratic collision operator \cite{Cai2015},
it can be expected that such a method can provide a quite accurate
approximation when the computational cost of the linearized collision
operator is acceptable. Another possiblity is to give each coefficient
a different decay rate:
\begin{displaymath}
\mQ^{**}_M[f](\bv) = P_M \mQ[P_M f](\bv) -
  \sum_{k_1 + k_2 + k_3 > M} \nu_M^{k_1 k_2 k_3}
    f_{k_1 k_2 k_3} H^{k_1 k_2 k_3}(\bv) \mM(\bv),
\end{displaymath}
and a possible choice of $\nu_M^{k_1 k_2 k_3}$ is the corresponding
term in the linearized collision operator:
\begin{displaymath}
\nu_M^{k_1 k_2 k_3} = -\frac{1}{k_1! k_2! k_3!} \int_{\mathbb{R}^3}
  \mathcal{L}[\varphi^{k_1 k_2 k_3}](\bv) H^{k_1 k_2 k_3}(\bv)
  \,\mathrm{d}\bv,
\end{displaymath}
where $\varphi^{k_1 k_2 k_3}(\bv) = H^{k_1 k_2 k_3}(\bv) \mM(\bv)$ is
the basis function. The effect of these finer approximations will be
studied in the future work.

By now, we have obtained a series of new collision models
\eqref{eq:approQ}. It can be expected that these models are better
approximations of the original quadratic operator than the simple
BGK-type models, especially when the non-equilibrium is strong and the
non-linearity takes effect. This will be observed in the numerical
examples.


\section{Numerical examples}
\label{sec:numerical}
In this section, we will show some results of our numerical
simulation. In all the numerical experiments, we adopt the newly
proposed collision operator \eqref{eq:approQ}, and solve the equation
\begin{displaymath}
\frac{\partial f}{\partial t} = \mQ_{M_0}[f]
\end{displaymath}
numerically for some positive integer $M_0$. This equation is solved
by the Galerkin spectral method with solution defined in the space
$\mF_M$, and $M$ is always chosen to be greater than $M_0$. For the
time discretization, we use the classical 4th-order Runge-Kutta method
in all the examples, and the time step is chosen as $\Delta t = 0.01$.

\subsection{BKW solution} \label{sec:BKW}
\begin{figure}[!ht]
\centering
\subfloat[Initial MDF $g(0,v_1)$\label{fig:ex1_init_1d}]{%
  \includegraphics[height=.25\textwidth]{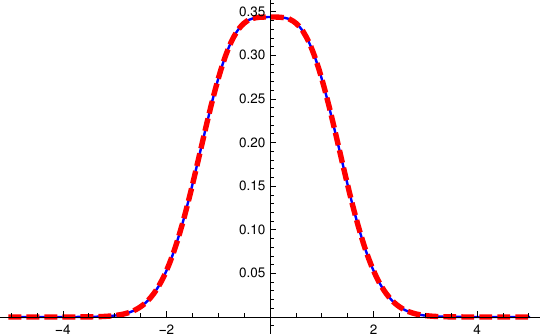}
}
\subfloat[Contours of $h(0,v_1,v_2)$\label{fig:ex1_init_2d_contour}]{%
  \includegraphics[height=.25\textwidth]{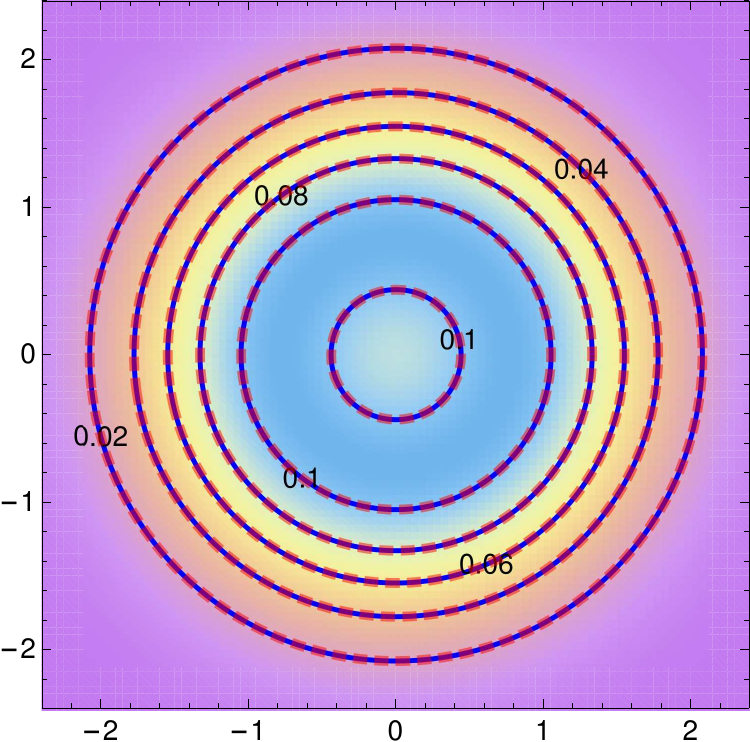}
}
\subfloat[Initial MDF $h(0,v_1,v_2)$\label{fig:ex1_init_2d}]{%
  \includegraphics[height=.25\textwidth]{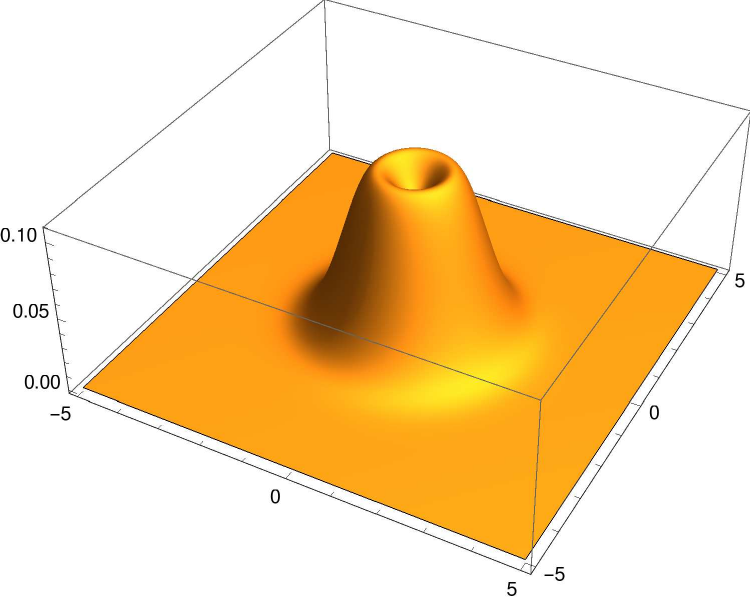}
}
\caption{Initial marginal distribution functions. In (a) and (b), the
  blue solid lines correspond to the exact solution, and the red
  dashed lines correspond to the numerical approximation. Figure (c)
  shows only the numerical approximation.}
\label{fig:ex1_init}
\end{figure}

For the Maxwell gas $\eta = 5$, the original spatially homogeneous
Boltzmann equation \eqref{eq:homo} admits an exact solution with
explicit expression:
\begin{displaymath}
f(t,\bv) = (2\pi \tau(t))^{-3/2} \exp \left( -\frac{|\bv|^2}{2\tau(t)} \right)
  \left[ 1 + \frac{1-\tau(t)}{\tau(t)}
    \left( \frac{|\bv|^2}{2\tau(t)} - \frac{3}{2} \right) \right],
\end{displaymath}
where
$\tau(t) = 1 - \exp\left(\frac{\pi}{3} \tilde{B}_2^5 (t+t_0)
\right)$.
In order that $f(\bv) \geqslant 0$ for all $t \in \mathbb{R}_+$ and
$\bv \in \mathbb{R}^3$, the parameter $t_0$ must satisfy
\begin{equation} \label{eq:t0}
-\frac{\pi}{3} \tilde{B}_2^5 t_0 \geqslant \log\left( \frac{5}{2} \right)
  \approx 0.916291.
\end{equation}
Here we choose $t_0$ such that the left hand side of \eqref{eq:t0}
equals to $0.92$. To ensure a good approximation of the initial
distribution function, we use $M = 20$ ($1771$ degrees of freedom) in
our simulation. For visualization purpose, we define the marginal
distribution functions (MDFs)
\begin{displaymath}
g(t,v_1) = \int_{\mathbb{R}} f(t,\bv) \,\mathrm{d}v_2 \,\mathrm{d}v_3, \qquad
h(t,v_1,v_2) = \int_{\mathbb{R}} f(t,\bv) \,\mathrm{d}v_3.
\end{displaymath}
The initial MDFs are plotted in Figure \ref{fig:ex1_init}, in which
the lines for exact functions and their numerical approximation are
hardly distinguishable.

\begin{figure}[!ht]
\centering
\subfloat[Profile of $g(0.2,v_1)$]{%
  \includegraphics[width=.3\textwidth]{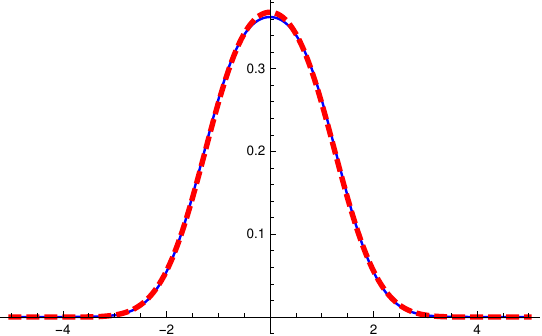}
} \hfill
\subfloat[Profile of $g(0.4,v_1)$]{%
  \includegraphics[width=.3\textwidth]{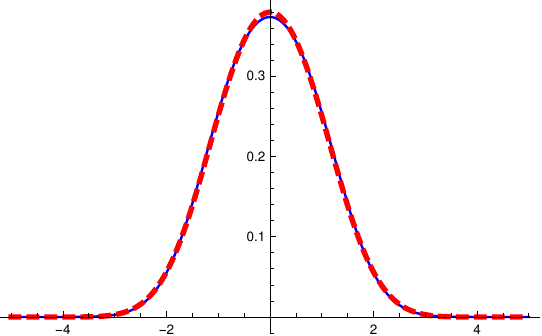}
} \hfill
\subfloat[Profile of $g(0.6,v_1)$]{%
  \includegraphics[width=.3\textwidth]{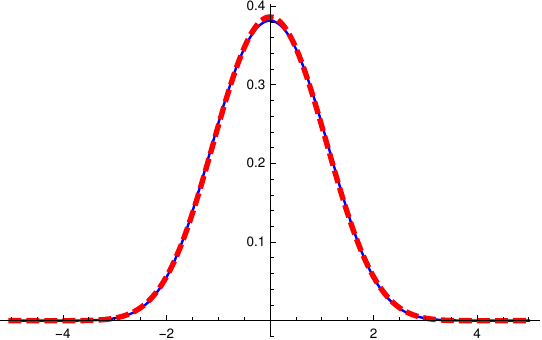}
} \\
\subfloat[Contours of $h(0.2,v_1,v_2)$]{%
  \includegraphics[width=.28\textwidth]{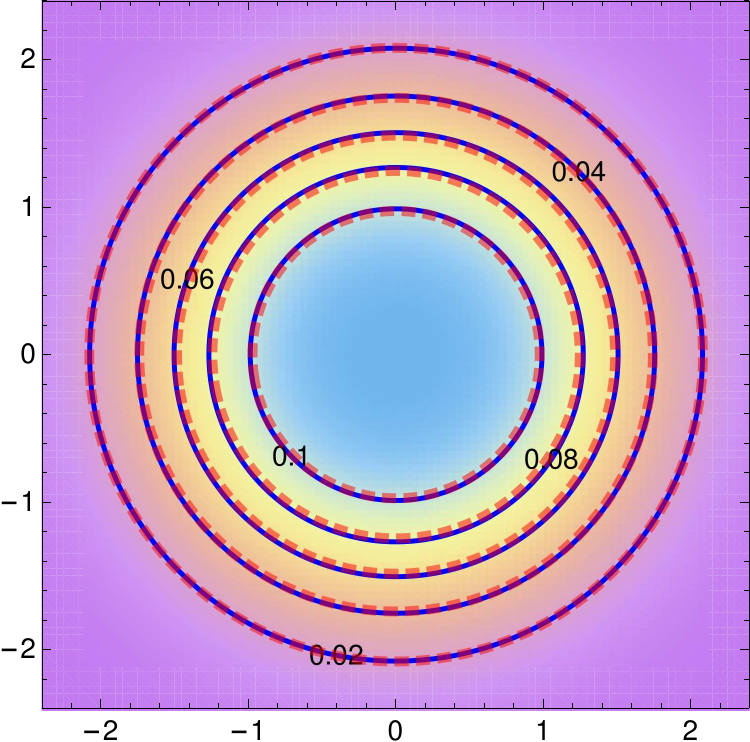}
} \hfill
\subfloat[Contours of $h(0.4,v_1,v_2)$]{%
  \includegraphics[width=.28\textwidth]{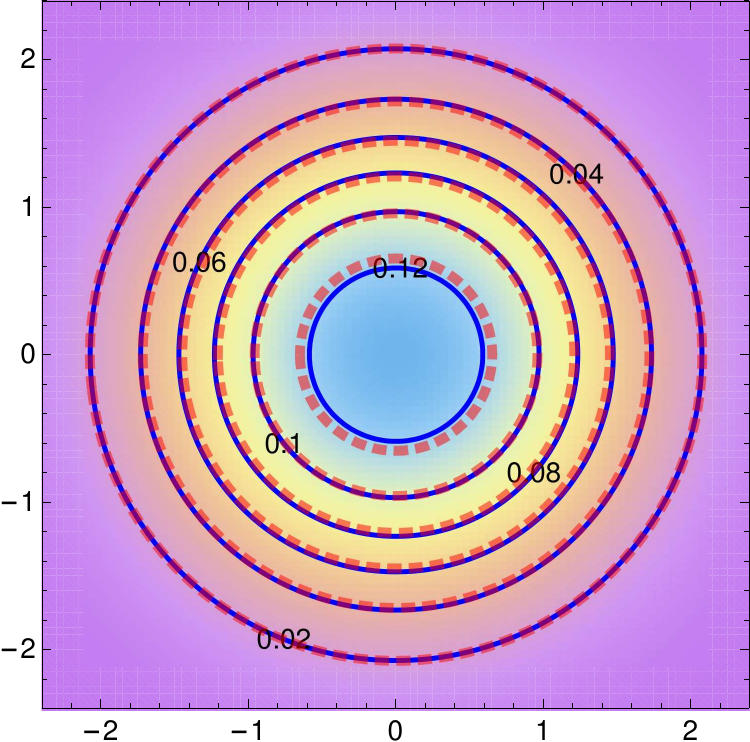}
} \hfill
\subfloat[Contours of $h(0.6,v_1,v_2)$]{%
  \includegraphics[width=.28\textwidth]{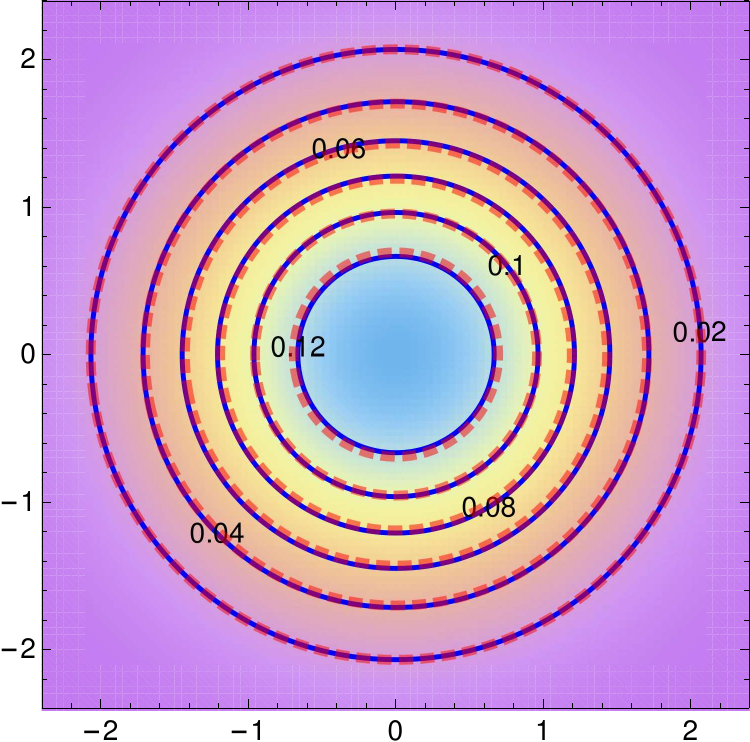}
}
\caption{Marginal distribution functions for $M_0 = 5$ at $t = 0.2$,
  $0.4$ and $0.6$. The blue lines correspond to the exact solution,
  and the red lines correspond to the numerical solutions.}
\label{fig:ex1_M0=5}
\end{figure}

\begin{figure}[!ht]
\centering
\subfloat[Profile of $g(0.2,v_1)$]{%
  \includegraphics[width=.3\textwidth]{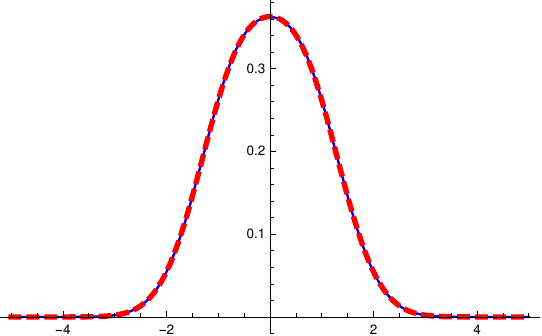}
} \hfill
\subfloat[Profile of $g(0.4,v_1)$]{%
  \includegraphics[width=.3\textwidth]{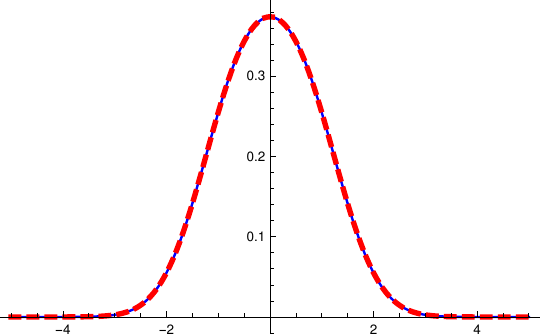}
} \hfill
\subfloat[Profile of $g(0.6,v_1)$]{%
  \includegraphics[width=.3\textwidth]{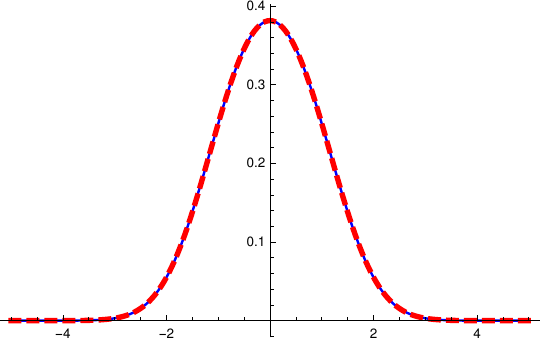}
} \\
\subfloat[Contours of $h(0.2,v_1,v_2)$]{%
  \includegraphics[width=.28\textwidth]{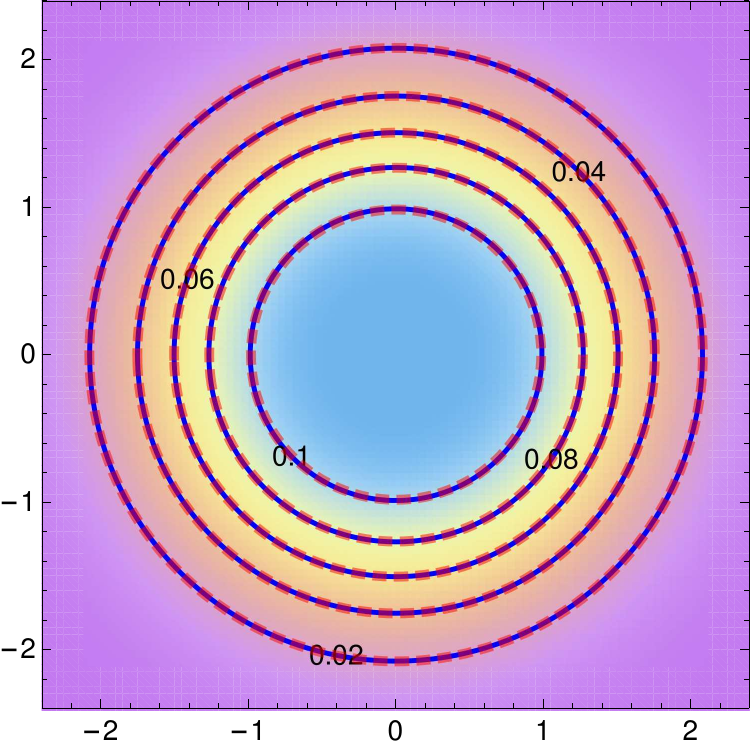}
} \hfill
\subfloat[Contours of $h(0.4,v_1,v_2)$]{%
  \includegraphics[width=.28\textwidth]{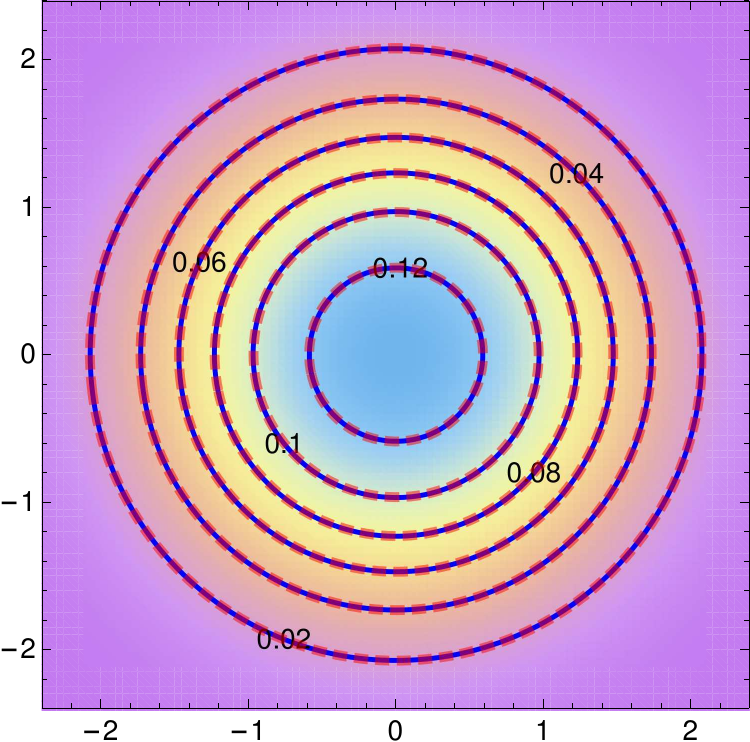}
} \hfill
\subfloat[Contours of $h(0.6,v_1,v_2)$]{%
  \includegraphics[width=.28\textwidth]{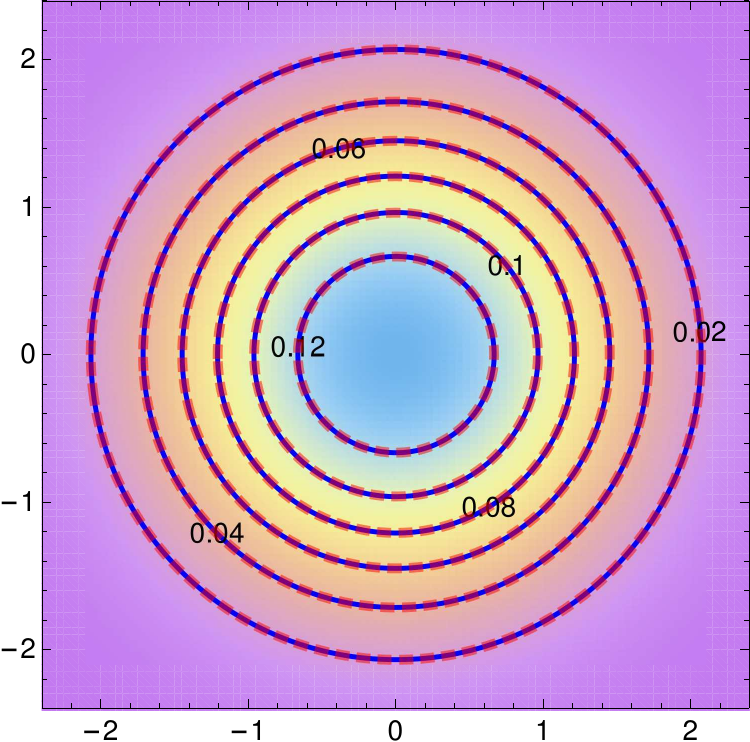}
}
\caption{Marginal distribution functions for $M_0 = 10$ at $t = 0.2$,
  $0.4$ and $0.6$. The blue lines correspond to the exact solution,
  and the red lines correspond to the numerical solutions.}
\label{fig:ex1_M0=10}
\end{figure}

\begin{figure}[!ht]
\centering
\includegraphics[width=.4\textwidth]{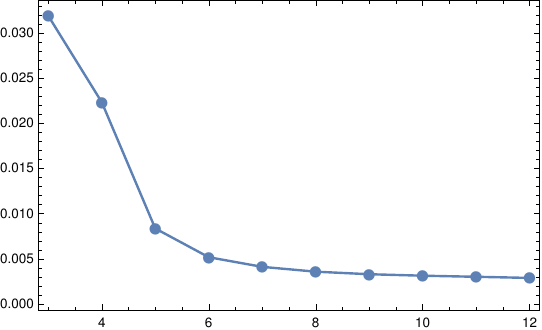}
\caption{The horizontal axis is the value of $M_0$, and the vertical
  axis is the value of $T_{M_0} / N_{M_0}^3$, where $T_{M_0}$ is the
  computational time (in milliseconds) for given $M_0$ and $N_{M_0}$
  is defined in \eqref{eq:coeN}.} 
\label{fig:ex1_comp_time}
\end{figure}

Numerical results for $t = 0.2$, $0.4$ and $0.6$ are given in Figures
\ref{fig:ex1_M0=5} and \ref{fig:ex1_M0=10}, respectively for $M_0 = 5$
and $M_0 = 10$. For $M_0 = 5$, the numerical solution provides a
reasonable approximation, but still with noticeable deviations, while
for $M_0 = 10$, the two solutions match perfectly in all cases. To
study the computational time, we run the simulation for
$M_0 = 3,\cdots,12$ until $t=5$ on a single CPU core with model Intel%
\textsuperscript{\textregistered} Core\texttrademark{} i7-7600U. The
relation between the computational time and the value of $M_0$ is
plotted in Figure \ref{fig:ex1_comp_time}. It can be seen that when
$M_0$ is large, the computational time is roughly proportional to the
cube of the number of degrees of freedom.  Note that the computational
time also includes the time for processing the coefficients of basis
functions with degree between $M_0 + 1$ and $M$. Although the time
complexity is only linear, when $M_0$ is small, the number of such
coefficients is quite large, and they have a significant contribution
to the total computational time. This explains why the curve in Figure
\ref{fig:ex1_comp_time} decreases fast for the first few points.
As a reference, we provide the average computational time for a
  single collsion operator in Table \ref{tab:comp_time}.

\begin{table}[!ht]
\centering
\begin{tabular}{ccccccccccc}
\hline
$M_0$ & 3 & 4 & 5 & 6 & 7 & 8 & 9 & 10 & 11 & 12 \\
\hline
Time (ms) & 0.128 & 0.479 & 0.734 & 1.535 &
  3.553 & 8.037 & 17.554 & 36.643 & 72.666 & 135.955 \\
\hline
\end{tabular}
\caption{Average computational time for a single collision operator
for different values of $M_0$.}
\label{tab:comp_time}
\end{table}

In Table \ref{tab:error}, we provide the $L^2$ and weighted $L^2$
error of the numerical solutions at $t=0.5$ and $t=1.0$. The notations
in the table are
\begin{displaymath}
E_M^{(1)} = \left( \int_{\mathbb{R}^3}
  |f_{\mathrm{num}}(\bv) - f_{\mathrm{exact}}(\bv)|^2 \,\mathrm{d}\bv
\right)^{1/2}, \quad E_M^{(2)} = \left( \int_{\mathbb{R}^3}
  |f_{\mathrm{num}}(\bv) - f_{\mathrm{exact}}(\bv)|^2 [\mM(\bv)]^{-1}
  \,\mathrm{d}\bv \right)^{1/2},
\end{displaymath}
where $f_{\mathrm{num}}$ is the numerical solution, and
$f_{\mathrm{exact}}$ is the exact solution. Four different choices of
$M_0$ ($M_0 = 5, 10, 15, 20$) and two different choices of $M$ ($M =
M_0$ and $M = 20$) are considered, from which we can see a rapid drop
of the numerical error as $M_0$ increases, indicating the spectral
accuracy. When $M_0 < 20$, the results for $M = 20$ are slightly more
accurate than the corresponding results for $M = M_0$, especially when
$M_0$ is small. We expect that such a property is useful when
simulating spatially inhomogeneous problems, for which the value of
$M_0$ cannot be too large due to the presence of the spatial
variables.

\begin{table}[!ht]
\centering
\renewcommand\arraystretch{1.2}
\footnotesize %
\begin{tabular}{|c|cccc|cccc|}
\hline
& \multicolumn{4}{c|}{$t=0.5$} & \multicolumn{4}{c|}{$t=1.0$} \\
\hline
$M_0$ & 5 & 10 & 15 & 20 & 5 & 10 & 15 & 20 \\
\hline
$E_{M_0}^{(1)}$ & $1.04 {\times} 10^{-2}$ & $5.40 {\times} 10^{-4}$ & $5.94 {\times} 10^{-5}$ & $1.90 {\times} 10^{-6}$ & $3.19 {\times} 10^{-3}$ & $6.09 {\times} 10^{-5}$ & $3.40 {\times} 10^{-6}$ & $3.89 {\times} 10^{-8}$ \\
$E_{M_0}^{(2)}$ & $7.46 {\times} 10^{-2}$ & $4.69 {\times} 10^{-3}$ & $5.57 {\times} 10^{-4}$ & $1.93 {\times} 10^{-5}$ & $2.52 {\times} 10^{-2}$ & $5.90 {\times} 10^{-4}$ & $3.50 {\times} 10^{-5}$ & $4.32 {\times} 10^{-7}$ \\
$E_{20}^{(1)}$ & $6.48 {\times} 10^{-3}$ & $3.71 {\times} 10^{-4}$ & $4.49 {\times} 10^{-5}$ & $1.90 {\times} 10^{-6}$ & $2.78 {\times} 10^{-3}$ & $5.53 {\times} 10^{-5}$ & $3.20 {\times} 10^{-6}$ & $3.89 {\times} 10^{-8}$ \\
$E_{20}^{(2)}$ & $5.05 {\times} 10^{-2}$ & $3.42 {\times} 10^{-3}$ & $4.31 {\times} 10^{-4}$ & $1.93 {\times} 10^{-5}$ & $2.28 {\times} 10^{-2}$ & $5.40 {\times} 10^{-4}$ & $3.31 {\times} 10^{-5}$ & $4.32 {\times} 10^{-7}$ \\
\hline
\end{tabular}
\caption{Numerical error for the BKW solution. $E_M^{(1)}$ is the
  $L^2$ error, and $E_M^{(2)}$ is the weighted $L^2$ error. See text
  for details.}
\label{tab:error}
\end{table}

\begin{figure}[!ht]
\centering
\subfloat[$f_{400}(t)$]{%
  \includegraphics[width=.4\textwidth]{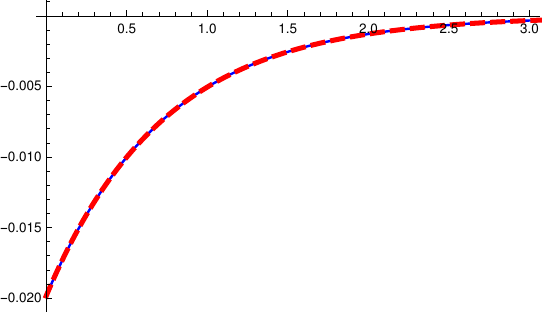}
} \hspace{20pt}
\subfloat[$f_{220}(t)$]{%
  \includegraphics[width=.4\textwidth]{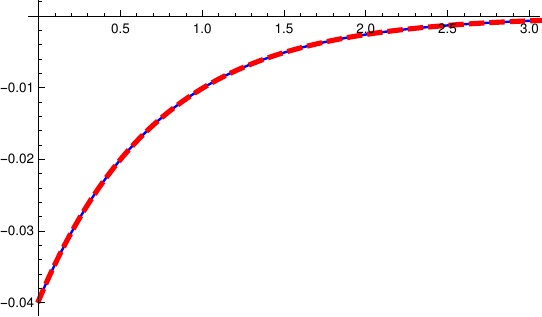}
}
\caption{The evolution of the coefficients. The blue lines correspond
  to the reference solution, and the red lines correspond to the
  numerical solution.}
\label{fig:ex1_moments}
\end{figure}

Now we consider the time evolution of the moments. By expanding the
exact solution into Hermite series, we get the exact solution for the
coefficients:
\begin{displaymath}
f_{k_1 k_2 k_3}(t) = \left\{\begin{array}{ll}
  \left[
    -\dfrac{1}{2} \exp \left( \dfrac{\pi}{3} \tilde{B}_2^{\eta} (t + t_0) \right)
  \right]^{\frac{k_1 + k_2 + k_3}{2}}
    \dfrac{1 - (k_1 + k_2 + k_3)/2}{(k_1/2)!(k_2/2)!(k_3/2)!}, &
  \text{if } k_1, k_2, k_3 \text{ are even}, \\[13pt]
  0, & \text{ otherwise}.
\end{array} \right.
\end{displaymath}
This exact solution can also be written in terms of initial
conditions as
\begin{displaymath}
f_{k_1 k_2 k_3}(t) = f_{k_1 k_2 k_3}(0) \exp \left(
  \frac{\pi}{6} \tilde{B}_2^{\eta} (k_1 + k_2 + k_3) t \right),
\end{displaymath}
from which one can clearly see that coefficients for higher-degree
polynomials decay faster, showing that a better modeling of the ``BGK
part'' may yield better results. Due to the symmetry of the
distribution function, the coefficients $f_{k_1 k_2 k_3}$ are zero for
any $t$ if $1 \leqslant k_1 + k_2 + k_3 \leqslant 3$. Hence we will
focus on the coefficients $f_{400}$ and $f_{220}$, which are the
fourth moments of the distribution function.  For Maxwell molecules,
the discrete kernel $A_{k_1 k_2 k_3}^{l_1 l_2 l_3 m_1 m_2 m_3}$ is
nonzero when $k_1 + k_2 + k_3 = l_1 + l_2 + l_3 + m_1 + m_2 + m_3$.
Therefore, for any $M \geqslant M_0 \geqslant 4$, the numerical
results for these two coefficients $f_{400}$ and $f_{220}$ are exactly
the same (regardless of round-off errors). Figure
\ref{fig:ex1_moments} gives the comparison between the numerical
solution and the exact solution for these two coefficients.  In both
plots, the two lines almost coincide with each other.

\subsection{Bi-Gaussian initial data}
\begin{figure}[!ht]
\centering
\subfloat[Initial MDF $g(0,v_1)$\label{fig:ex2_init_1d}]{%
  \includegraphics[height=.25\textwidth]{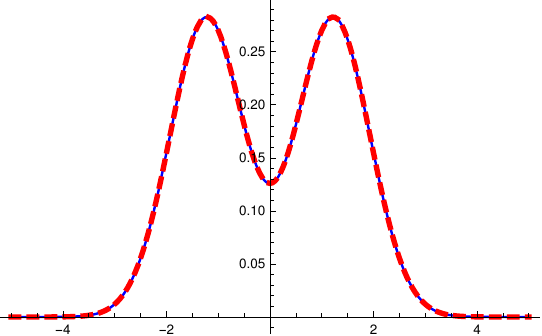}
}
\subfloat[Contours of $h(0,v_1,v_2)$\label{fig:ex2_init_2d_contour}]{%
  \includegraphics[height=.25\textwidth]{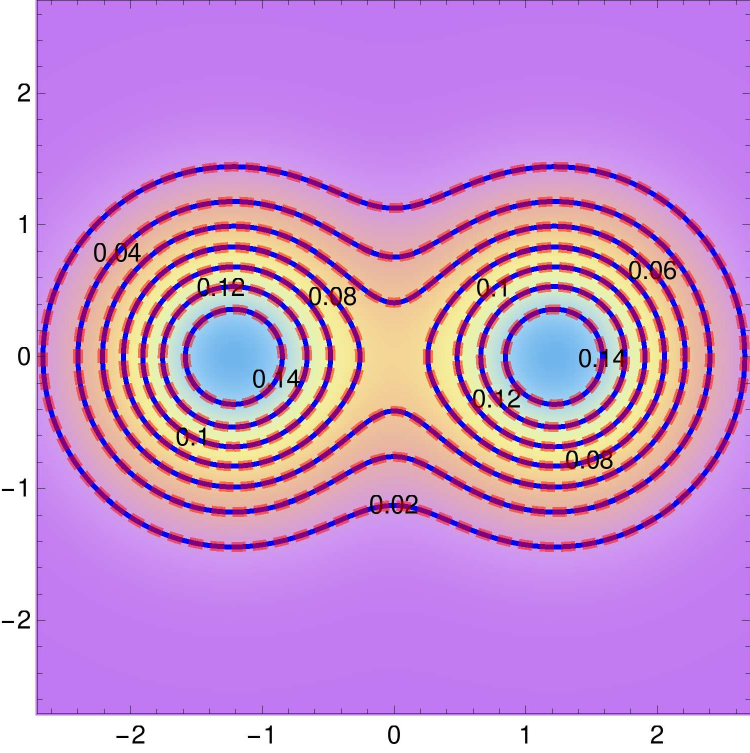}
}
\subfloat[Initial MDF $h(0,v_1,v_2)$\label{fig:ex2_init_2d}]{%
  \includegraphics[height=.25\textwidth]{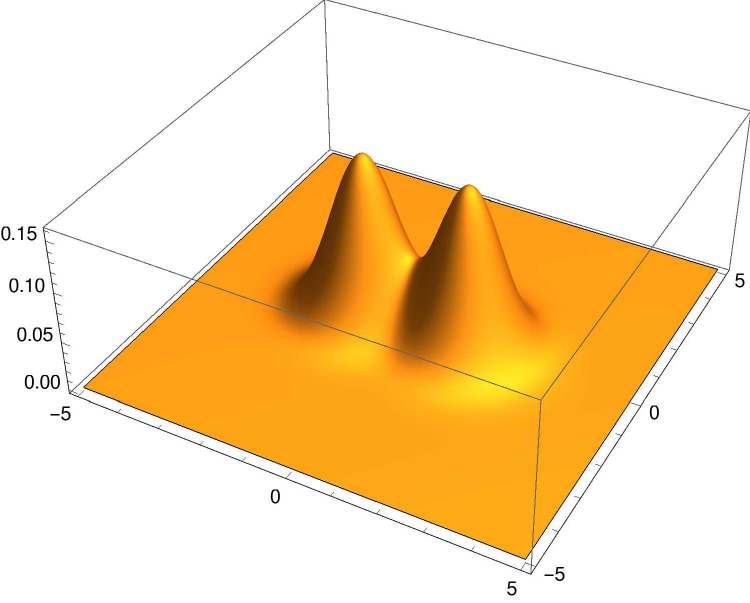}
}
\caption{Initial marginal distribution functions. In (a) and (b), the
  blue solid lines correspond to the exact solution, and the red
  dashed lines correspond to the numerical approximation. Figure (c)
  shows only the numerical approximation.}
\label{fig:ex2_init}
\end{figure}

In this example, we perform the numerical test for hard potential
$\eta = 10$. The initial distribution function is
\begin{displaymath}
f(0,\bv) = \frac{1}{2\pi^{3/2}} \left[
  \exp \Big( -(v_1 + \sqrt{3/2})^2 + v_2^2 + v_3^2 \Big) +
  \exp \Big( -(v_1 - \sqrt{3/2})^2 + v_2^2 + v_3^2 \Big)
\right].
\end{displaymath}
Again, in all our numerical tests, we use $M = 20$ which gives a good
approximation of the initial distribution function (see Figure
\ref{fig:ex2_init}).

\begin{figure}[!ht]
\centering
\subfloat[$t=0.3$]{%
  \includegraphics[width=.33\textwidth]{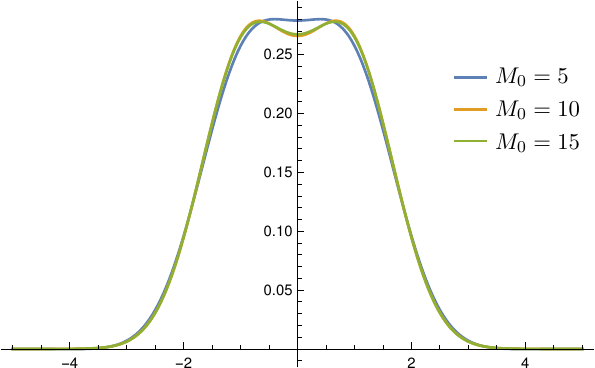}
}
\subfloat[$t=0.6$]{%
  \includegraphics[width=.33\textwidth]{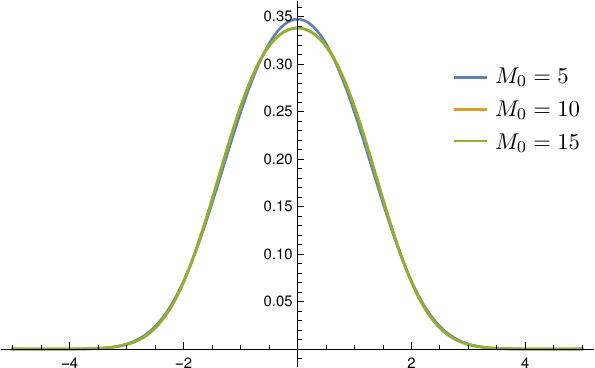}
}
\subfloat[$t=0.9$]{%
  \includegraphics[width=.33\textwidth]{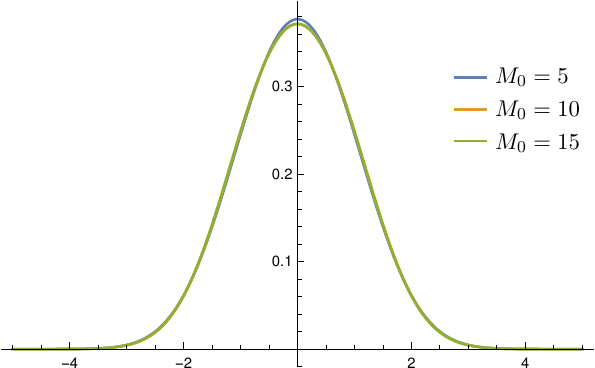}
}
\caption{Marginal distribution functions at different times.}
\label{fig:ex2_1d}
\end{figure}

\begin{figure}[!ht]
\centering
\subfloat[$t=0.3$]{%
  \includegraphics[width=.3\textwidth]{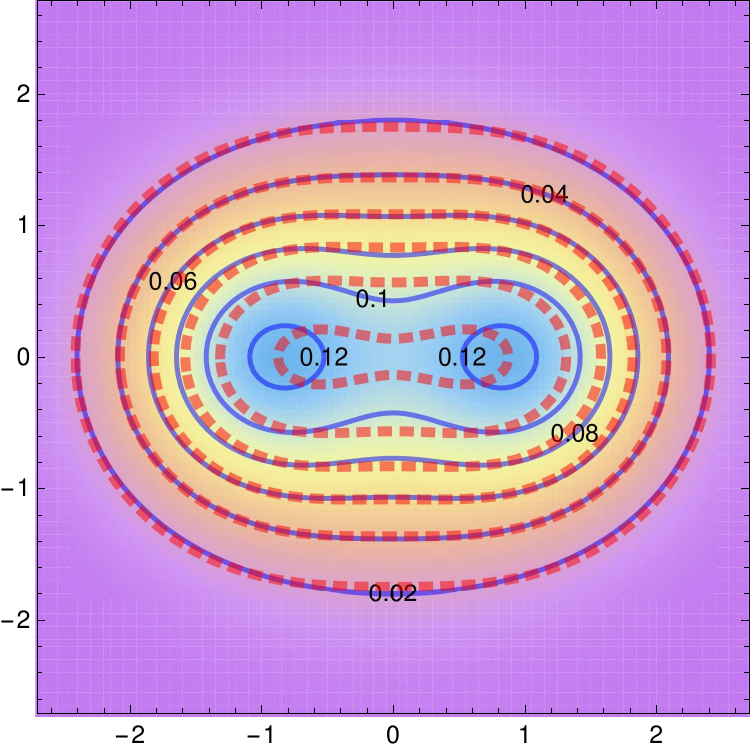}
} \hfill
\subfloat[$t=0.6$]{%
  \includegraphics[width=.3\textwidth]{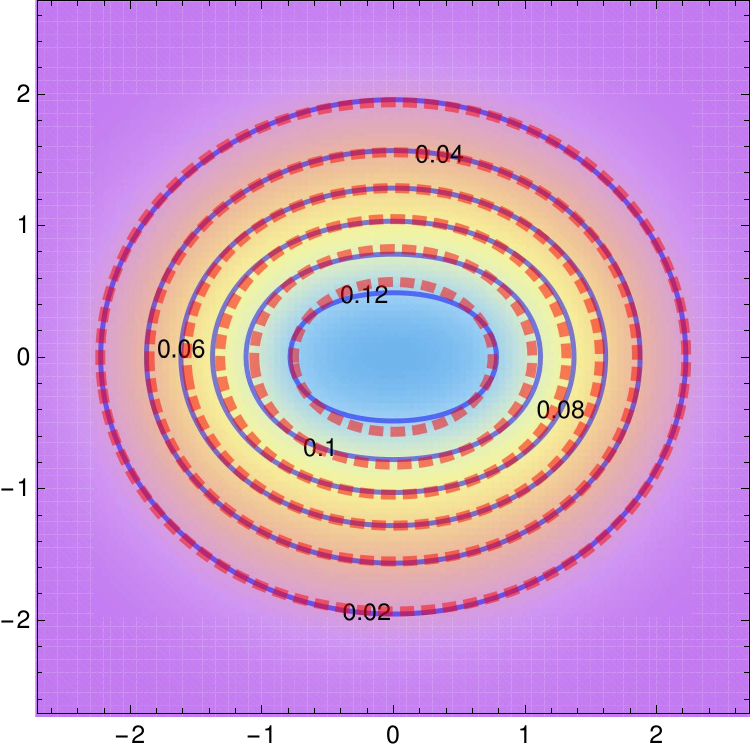}
} \hfill
\subfloat[$t=0.9$]{%
  \includegraphics[width=.3\textwidth]{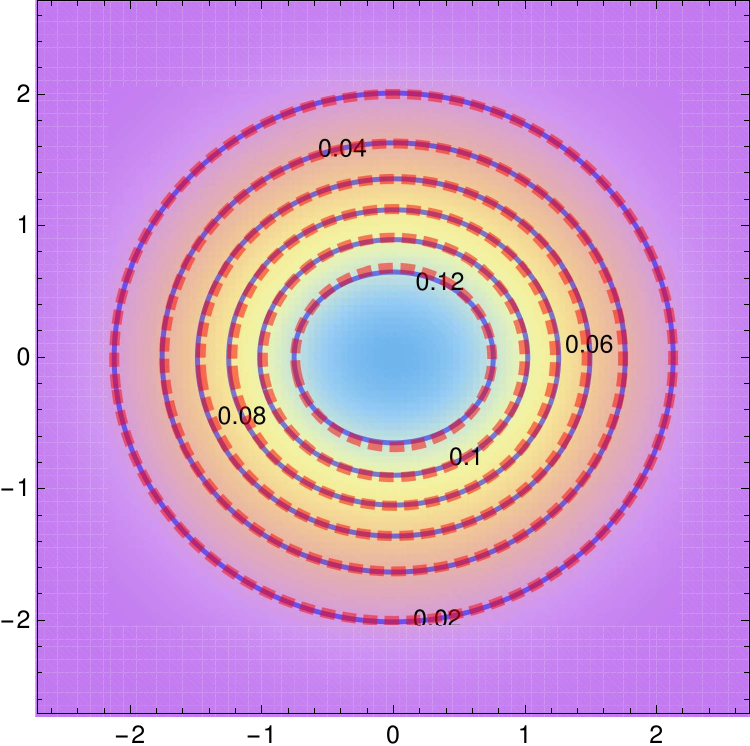}
}
\caption{Comparison of numerical results using $M_0 = 5$ and
  $M_0 = 15$. The blue contours and the red dashed contours are
  respectively the results for $M_0 = 5$ and $M_0 = 15$.}
\label{fig:ex2_2d_M0=5_15}
\end{figure}

\begin{figure}[!ht]
\centering
\subfloat[$t=0.3$]{%
  \includegraphics[width=.3\textwidth]{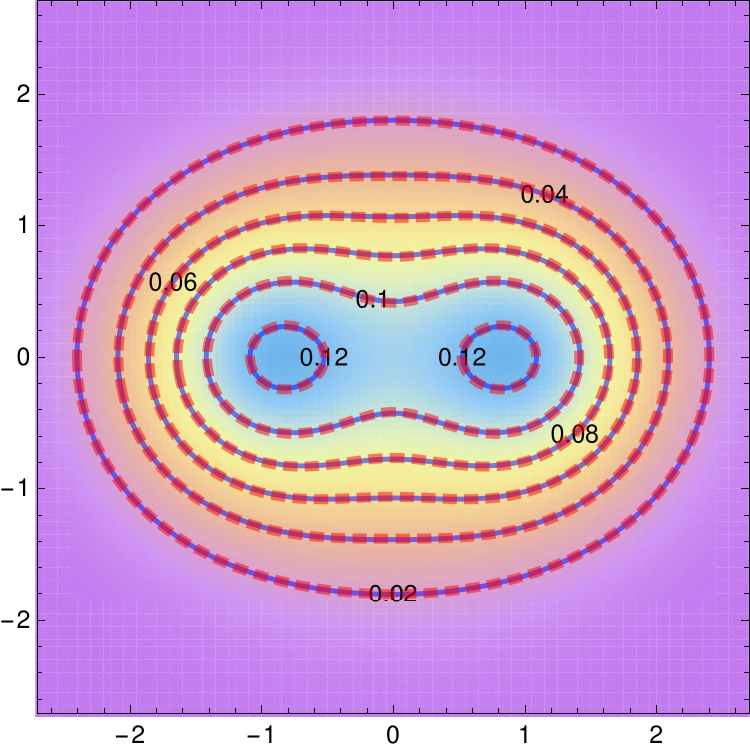}
} \hfill
\subfloat[$t=0.6$]{%
  \includegraphics[width=.3\textwidth]{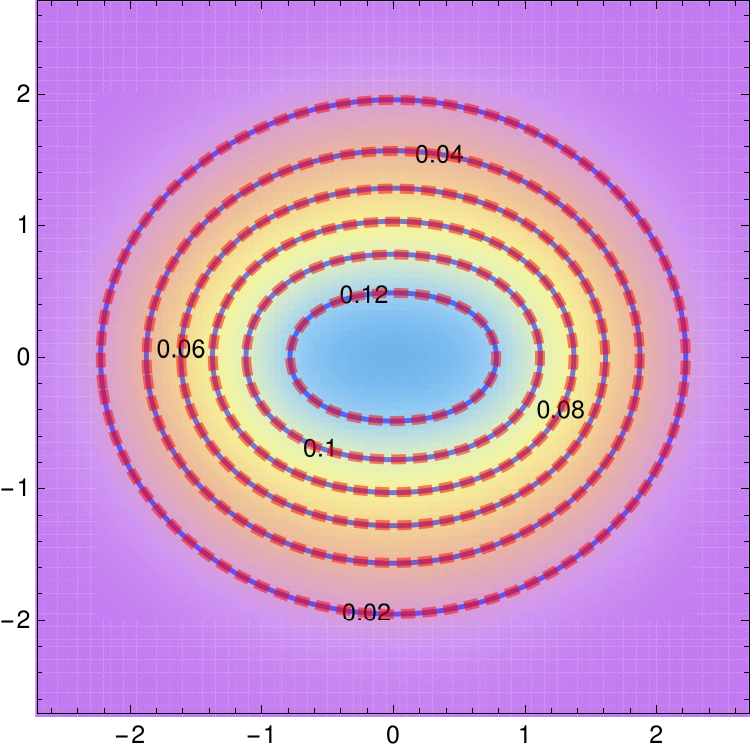}
} \hfill
\subfloat[$t=0.9$]{%
  \includegraphics[width=.3\textwidth]{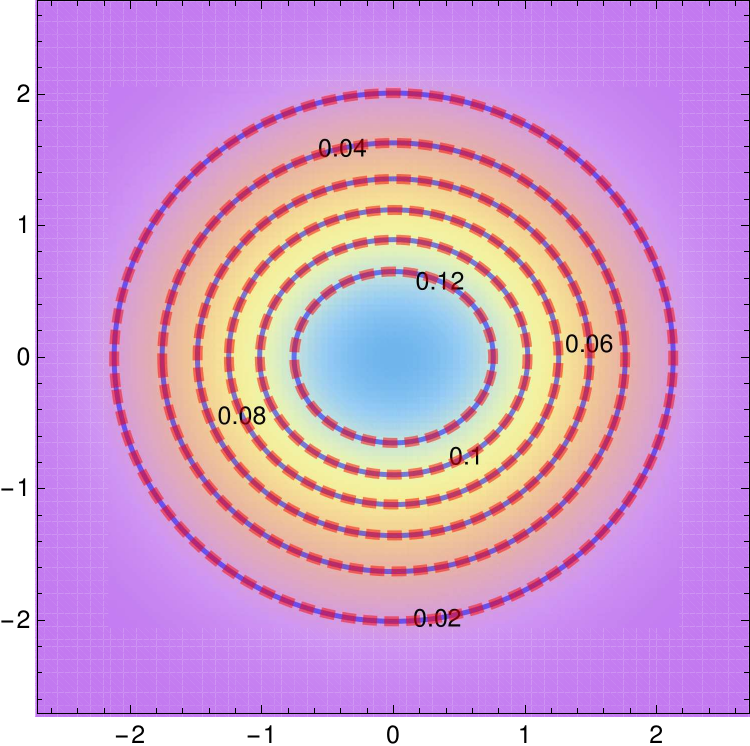}
}
\caption{Comparison of numerical results using $M_0 = 10$ and
  $M_0 = 15$. The blue contours and the red dashed contours are
  respectively the results for $M_0 = 10$ and $M_0 = 15$.}
\label{fig:ex2_2d_M0=10_15}
\end{figure}

For this example, we consider the three cases $M_0 = 5, 10, 15$, and
the corresponding one-dimensional marginal distribution functions at
$t = 0.3$, $0.6$ and $0.9$ are given in Figure \ref{fig:ex2_1d}. In
all the results, the lines for $M_0 = 10$ and $M_0 = 15$ are very
close to each other. Due to the fast convergence of the spectral
method, it is believable that $M_0 = 10$ can already provide a very
good approximation. To get a clearer picture, similar comparison of
two-dimensional results are also provided in Figure
\ref{fig:ex2_2d_M0=5_15} and \ref{fig:ex2_2d_M0=10_15}.

\begin{figure}[!ht]
\centering
\includegraphics[width=.4\textwidth]{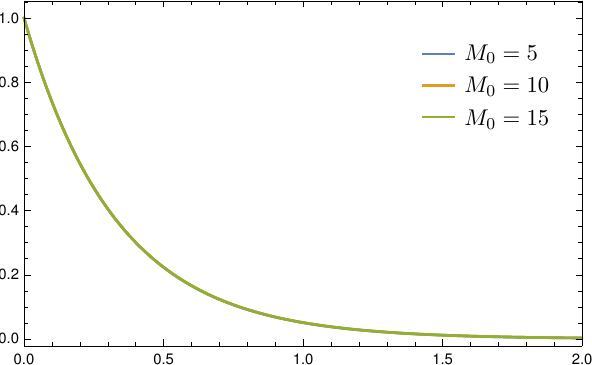}
\caption{Evolution of $\sigma_{11}(t)$. Three lines are on top of each other.}
\label{fig:ex2_sigma11}
\end{figure}

Now we consider the evolution of the moments. In this example, we
always have $\sigma_{11} = -2\sigma_{22} = -2\sigma_{33}$ and
$q_1 = q_2 = q_3 = 0$.  Therefore we focus only on the evolution of
$\sigma_{11}$, which is plotted in Figure \ref{fig:ex2_sigma11}. It
can be seen that three tests give almost identical results. Even for
$M_0 = 5$, while the distribution function is not approximated very
well, the evolution of the stress tensor is almost exact.

\subsection{Discontinuous initial data}
\begin{figure}[!ht]
\centering
\begin{tabular}{rr}
\subfloat[$\sigma_{11}(t)$ ($\eta = 10$)]{%
  \includegraphics[scale=.6]{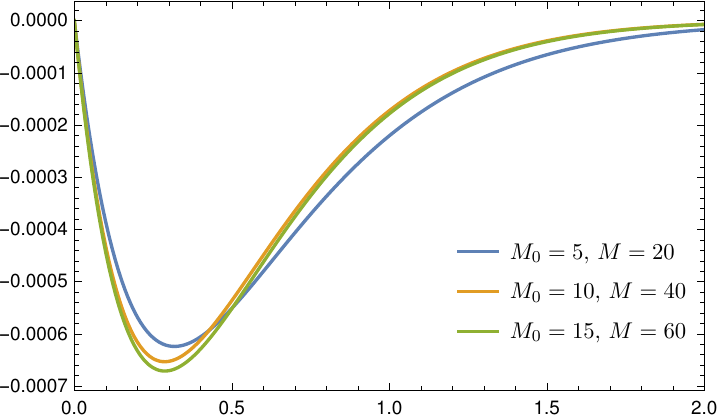}
} &
\subfloat[$\sigma_{11}(t)$ ($\eta = 3.1$)]{%
  \includegraphics[scale=.6]{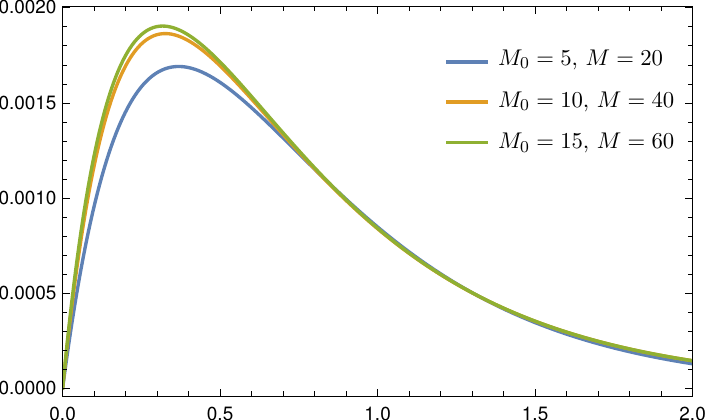}
} \\[20pt]
\subfloat[$\sigma_{22}(t)$ ($\eta = 10$)]{%
  \includegraphics[scale=.6]{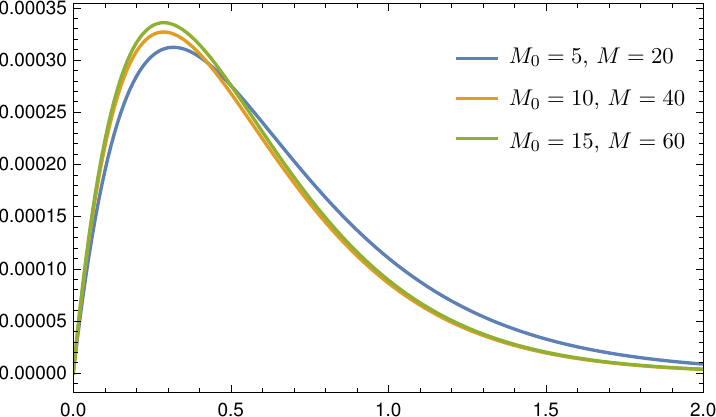}
} &
\subfloat[$\sigma_{22}(t)$ ($\eta = 3.1$)]{%
  \includegraphics[scale=.6]{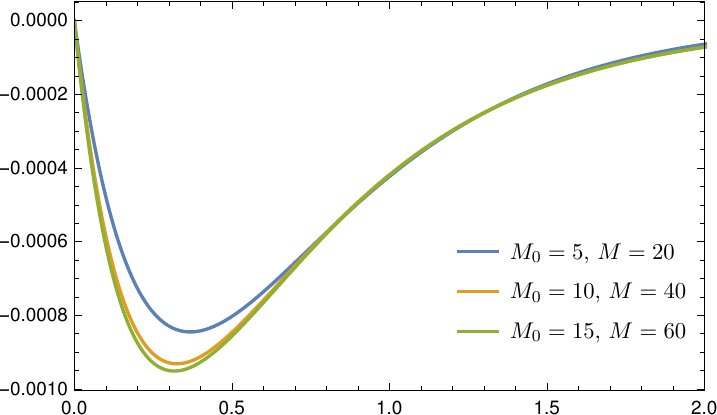}
} \\[20pt]
\subfloat[$q_1(t)$ ($\eta = 10$)]{%
  \includegraphics[scale=.6]{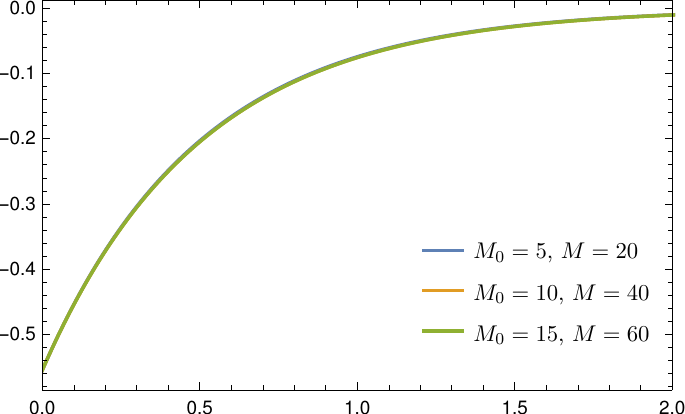}
} &
\subfloat[$q_1(t)$ ($\eta = 3.1$)]{%
  \includegraphics[scale=.6]{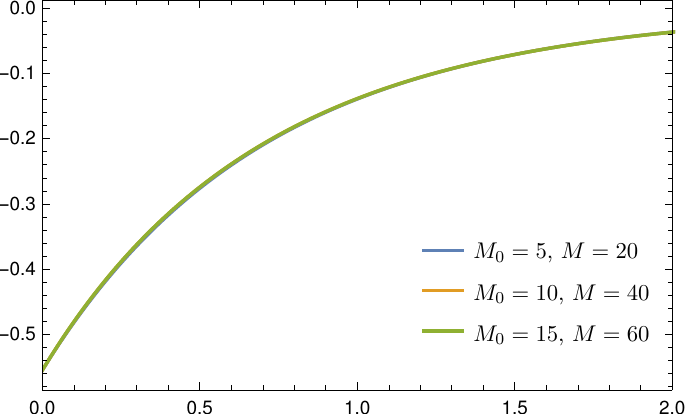}
}
\end{tabular}
\caption{Evolution of the stress and the heat flux. The left column
  shows the results for $\eta = 10$, and the right column shows the
  results for $\eta = 3.1$. In the right column, the horizontal axes
  are the scaled time (see \eqref{eq:scaling} and the context for
  details).}
\label{fig:ex3_moments}
\end{figure}

Here we consider the problem with a discontinuous initial condition:
\begin{displaymath}
f(0,\bv) = \left\{ \begin{array}{ll}
  \dfrac{\sqrt[4]{2} (2-\sqrt{2})}{\pi^{3/2}}
    \exp \left( -\dfrac{|\bv|^2}{\sqrt{2}} \right), & \text{if } v_1 > 0, \\[10pt]
  \dfrac{\sqrt[4]{2} (2-\sqrt{2})}{4\pi^{3/2}}
    \exp \left( -\dfrac{|\bv|^2}{2\sqrt{2}} \right), & \text{if } v_1 < 0.
\end{array} \right.
\end{displaymath}
We refer the readers to \cite{Cai2015} for the graphical profile of
this initial value. As a spectral method, the truncated expansion
\eqref{eq:f_h} is difficult to capture an accurate profile of a
discontinuous function. Therefore, we focus only on the evolution of
the moments. The left column of Figure \ref{fig:ex3_moments} shows the
numerical results for $\eta = 10$ with different choices of $M_0$ and
$M$. All the numerical tests show that the magnitude of the stress
components $\sigma_{11}$ and $\sigma_{22}$, which are initially zero,
increases to a certain number before decreasing again. Such phenomenon
cannot be captured by the simple BGK-type models. The lines
corresponding to the results of $M_0 = 10$, $M = 40$ and $M_0 = 15$,
$M = 60$ are very close to each other, which indicates that they might
be very close to the exact solution. For the case $M_0 = 5$, $M = 20$,
although an obvious error can be observed, the trends of the evolution
are qualitatively correct, and thus the corresponding collision model
$\mQ_5[f]$ may also be used as a better alternative to the BGK-type
models. For the heat flux $q_1$, the three results are hardly
distinguishable.

The right column of Figure \ref{fig:ex3_moments} gives the same
moments for the soft potential $\eta = 3.1$. For comparison purpose,
the horizontal axes are the scaled time $t_s = t / \tau$, where
\begin{equation} 
  \label{eq:scaling}
  \tau = \frac{4^{\frac{2}{\eta-1} - \frac{2}{9}} \tilde{B}_2^{10} \Gamma(34/9)}
  {\tilde{B}_2^{\eta} \Gamma(4-2/(\eta-1))} \approx 2.03942.
\end{equation}
By such scaling, the two models $\eta = 10$ and $\eta = 3.1$ have the
same mean relaxation time near equilibrium. The two columns in Figure
\ref{fig:ex3_moments} show quite different behavior for different
collision models, while both numerical results indicate the high
efficiency of this method in capturing the behavior of the moments.

\section{Concluding remarks and comparison with similar works}
\label{sec:conclusion} 
This work aims at an affordable way to model and simulate the binary
collision between gas molecules. Our new attempt is an intermediate
approach between a direct discretization of the quadratic Boltzmann
collision operator and simple modelling methods like BGK-type
operators. In detail, we first focus on the relatively important
physical quantities, which are essentially the first few coefficients
in the Hermite expansion, and use an intricate and accurate way to
describe their evolution. The strategy comes from the discretization
of the quadratic collision operator. For the less important
quantities, we borrow the idea of the BGK-type operators and let them
converge to the equilibrium at a constant rate.  Although the first
part is computationally expensive, we can restrict the number of
degrees of freedom such that the computational cost is acceptable. The
accuracy of such a model depends apparently on the size of the
accurately modelled part.

In the literature, there are already some works implementing
the Hermite spectral method using different algorithms, among which
\cite{Gamba2018, Kitzler2019} is essentially the same as ours. The
difference is the implementation: the work \cite{Gamba2018} uses
orthogonal polynomials based on spherical coordinates in the
three-dimensional Euclidean space, while we use the orthogonal
polynomials based on the Cartesian coordinates; the work
\cite{Kitzler2019} uses the same orthogonal polynomials as ours, while
the proposed computational cost in \cite{Kitzler2019} is $O(M^7)$.
Compared with \cite{Gamba2018}, in which the coefficients are computed
numerically, we can compute all the coefficients $A_{k_1 k_2 k_3}^{i_1
i_2 i_3, j_1 j_2 j_3}$ almost exactly, except for the one-dimensional
integration in \eqref{eq:tB_exact}. Compared with the algorithm in
\cite{Kitzler2019}, our method has a higher time complexity $O(M^9)$
if the full quadratic collision operator is used. Despite this, one
can directly compare the computational time for both algorithms. It
seems that our algorithm is still faster when $M$ is small, due to a
larger constant hidden in front of their computational cost $M^7$. One
obvious deficiency of our algorithm is the memory cost as listed in
Table \ref{tab:memory}. We need $O(M^9)$ while \cite{Kitzler2019}
needs only $O(M^4)$. The reason of such a difference is that
the work \cite{Kitzler2019} has shifted most of our calculation in the
appendix to the online computation, whereas we store these
intermediate results in memory. This leads to different memory cost
for the two algorithms. Moving these computations online also makes it
possible to reduce the time complexity. Thus our algorithm to compute
the full quadrature collision operator will eventually be slower as
$M$ increases. Therefore in Section \ref{sec:convergence}, we proposed
a remedy to allow computations with a large $M$.

Another related work is \cite{Grohs2017}, where the basis functions
are chosen such that the discretization is in the $L^2$ space instead
of the weighted $L^2$ space. One advantage of this method is that
$L^2(\mathbb{R}^3)$ is a large space, and more distribution functions
can be included to the framework. However, since the coefficients in
the expansion are not directly related to the moments, and the trick
of cost reduction in Section \ref{sec:convergence} is not applicable.

Our numerical examples show that our method can efficiently capture
the evolution of lower-order moments in the spatially homogeneous
Boltzmann equation. The method should be further validated in the
numerical tests for the full Boltzmann equation with spatial
variables, by which one can probably get a proper {\it a priori}
estimation of $M_0$. Some preliminary applications to several
benchmark problems have been done in \cite{Hu2018}, and more
experiments are to be carried out in future works. Besides, we are also working on a
better choice of the ``BGK part'' in our collision model and the
reduction of the computational cost for the quadratic part.

\section*{Acknowledgements}
We would like to thank Prof. Manuel Torrilhon at RWTH Aachen
University, Germany for motivating this research project and Prof. Ruo
Li at Peking University, China for the valuable suggestions to this
research project. Zhenning Cai is supported by National University of
Singapore Startup Fund under Grant No. R-146-000-241-133. Yanli Wang
is supported by the National Natural Scientific Foundation of China
(Grant No.  11501042 and 91630310).

\appendix
\section{Proof of Theorem \ref{thm:coeA}}
\label{sec:thm_proof}

In order to prove Theorem \ref{thm:coeA}, we first introduce the lemma
below:
\begin{lemma}
\label{thm:Hermit_cv}
Let $\bv = \bh + \bg/2$ and $\bw = \bh - \bg/2$. It holds that
\begin{displaymath}
\begin{split}
& H^{k_1 k_2 k_3}(\bv) H^{l_1 l_2 l_3}(\bw) = \\
& \quad \sum_{k_1'+l_1' = k_1+l_1}
  \sum_{k_2'+l_2' = k_2+l_2} \sum_{k_3'+l_3' = k_3+l_3}
  a_{k_1'l_1'}^{k_1 l_1} a_{k_2'l_2'}^{k_2 l_2} a_{k_3'l_3'}^{k_3 l_3}
  H^{k_1' k_2' k_3'}(\sqrt{2}\bh)
  H^{l_1' l_2' l_3'} \left( \frac{\bg}{\sqrt{2}} \right),
\end{split}
\end{displaymath}
where the coefficients $a_{k_s'l_s'}^{k_s l_s}$, $s = 1,2,3$ are
defined in \eqref{eq:coea}.
\end{lemma}

{\renewcommand\proofname{Proof of Lemma \ref{thm:Hermit_cv}}
  \begin{proof}
    First, it is easy to verify that
    $ \exp\left(-\frac{|\bv|^2 + |\bw|^2}{2}\right) =
    \exp\left(-\left(|\bh|^2 + \frac{|\bg|^2}{4}\right)\right)$
    and $\dd\bv\dd \bw = \dd \bg \dd \bh$. Based on the orthogonality
    of the Hermite polynomials \eqref{eq:Her_orth}, we just need to
    prove
    \begin{equation} \label{eq:zeta_a}
    \czeta = \left\{ \begin{array}{ll}
      k_1'! k_2'! k_3'! l_1'! l_2'! l_3'!
      a_{k_1'l_1'}^{k_1 l_1} a_{k_2'l_2'}^{k_2 l_2}a_{k_3'l_3'}^{k_3 l_3},
        & \text{if } k_s + l_s = k_s' + l_s', \quad \forall s=1,2,3,
	\\
      0, & \text{otherwise},
    \end{array} \right.
    \end{equation}
    where the left hand side is defined as
    \begin{equation}
      \label{eq:coeH}
      \begin{aligned}
        \czeta := & \int_{\bbR^3}\int_{\bbR^3}
        H^{k_1 k_2 k_3}(\bv) H^{l_1 l_2 l_3}(\bw)
        H^{k_1' k_2' k_3'}(\sqrt{2}\bh) H^{l_1' l_2' l_3'} \left(
          \frac{\bg}{\sqrt{2}} \right) \exp\left(-\frac{|\bv|^2 +
            |\bw|^2}{2}\right) \dd \bv \dd \bw.
      \end{aligned}
    \end{equation}
    By the general Leibniz rule, we have the following relation for
    the derivatives of with respect to $\bv, \bw$ and $\bg, \bh$:
    \begin{equation}
      \label{eq:bin_diff}
      \frac{\partial^{k_s + l_s}}{\partial v_s^{k_s}\partial w_s^{l_s}} = 
      \sum_{i_s =0}^{k_s}\sum_{j_s=0}^{l_s}\binom{k_s}{i_s}\binom{l_s}{j_s}
      \frac{(-1)^{l_s - j_s}}{2^{i_s+j_s}}
      \frac{\partial^{k_s}}{\partial
        h_s^{i_s+j_s}}\frac{\partial^{l_s}}{\partial g_s^{i'_s+j'_s}}, \qquad
      i'_s = k_s - i_s, \quad j'_s = l_s - j_s, \quad s   = 1, 2, 3.
    \end{equation}
    Then, following the definition of Hermite polynomials
    \eqref{eq:basis} and \eqref{eq:bin_diff}, and using integration
    by parts, we arrive at
    \begin{equation}
      \label{eq:coeH_detail}
      \begin{aligned}
        & \czeta = \int_{\bbR^3}\int_{\bbR^3}
        \exp\left(-\left(|\bh|^2+\frac{|\bg|^2}{4}\right)\right) \times \\
        & \qquad \prod_{s=1}^3\left(\sum_{i_s
            =0}^{k_s}\sum_{j_s=0}^{l_s}\binom{k_s}{i_s}\binom{l_s}{j_s}
          \frac{(-1)^{l_s - j_s}}{2^{i_s+j_s}}\frac{ \partial^{k_s + l_s}
          }{\partial h_s^{i_s+j_s}\partial g_s^{i'_s+j'_s}} \right)H^{k_1'
          k_2' k_3'}(\sqrt{2}\bh)H^{l_1' l_2' l_3'}
        \left(\frac{\bg}{\sqrt{2}} \right) \dd \bh \dd \bg.
      \end{aligned}
    \end{equation}
    From the orthogonality of Hermite polynomials and the
    differentiation relation
    \begin{equation} \label{eq:Hermite_diff}
    \frac{\partial}{\partial v_s} H^{k_1 k_2 k_3}(\bv) = \left\{
      \begin{array}{ll}
      0, & \text{if } k_s = 0, \\
      k_s H^{k_1-\delta_{1s}, k_2-\delta_{2s}, k_3-\delta_{3s}}(\bv),
        & \text{if } k_s > 0,
      \end{array}
    \right.
    \end{equation}
    it holds that \eqref{eq:coeH_detail} is nonzero only when $i_s +
    j_s = k'_s$, $i'_s + j'_s = l'_s$, $s = 1,2,3$, which means
    \begin{equation}
      \label{eq:index_H}
      k_s + l_s = k'_s + l'_s, \qquad \forall s = 1,2,3.
    \end{equation}
    When \eqref{eq:index_H} holds, we can apply
    \eqref{eq:Hermite_diff} to \eqref{eq:coeH_detail} and get
      \begin{equation}
        \label{eq:coeH_detail1}
        \begin{aligned}
          \frac{\czeta}{
            k'_1!k'_2!k'_3!l'_1!l'_2!l'_3!}  & =
          \prod_{s=1}^3\left(\sum_{i_s =0}^{k_s}\sum\limits_{j_s=0,
              i_s + j_s = k'_s}^{l_s}\binom{k_s}{i_s}\binom{l_s}{j_s}
            \frac{(-1)^{l'_s - k_s + i_s}}{2^{k'_s}} 2^{\frac{k'_s -
                l'_s}{2}}\right) = a_{k_1'l_1'}^{k_1 l_1}
          a_{k_2'l_2'}^{k_2 l_2} a_{k_3'l_3'}^{k_3 l_3}.
        \end{aligned}
      \end{equation}
    Thus \eqref{eq:zeta_a} is shown, which completes the proof of the
    lemma.
  \end{proof}
}

\begin{corollary}
  \label{thm:Hermit_split}
  Let $\bv = \bh + \bg/2$. We have
  \begin{displaymath}
    H^{k_1 k_2 k_3}(\bv) = \sum_{l_1 + m_1 = k_1}
    \sum_{l_2 + m_2 = k_2} \sum_{l_3 + m_3 = k_3}
    \frac{2^{-k/2} k_1! k_2! k_3!}{l_1! l_2! l_3! m_1! m_2! m_3!}
    H^{l_1 l_2 l_3}(\sqrt{2} \bh)
    H^{m_1 m_2 m_3} \left( \frac{\bg}{\sqrt{2}} \right).
  \end{displaymath}
\end{corollary}
{\renewcommand\proofname{Proof of Corollary \ref{thm:Hermit_split}} 
  \begin{proof}
  This corollary is just a special case of Lemma \ref{thm:Hermit_cv}
  when $l_1 = l_2 = l_3 = 0$.
  \end{proof}
}

{\renewcommand\proofname{Proof of Theorem \ref{thm:coeA}}
  \begin{proof}
    Let $\bw = \bv', \bw_1 = \bv'_1, \bs = \bw - \bw_1$ and define the
    unit vector $\tilde{\bn}$ as $\tilde{\bn} = -(\bg\sin\chi /|\bg| +
    \bn\cos\chi)$. It holds that
    \begin{equation}
      \begin{gathered}
        \label{eq:conser_v}
        |\bv|^2 + |\bv_1|^2 = |\bw|^2 +
        |\bw_1|^2, \qquad \dd \bv\dd \bv_1 = \dd \bw \dd \bw_1, \qquad
        |\bs|= |\bg|, \qquad \bs \cdot \bn_{\bw} = 0, \\
        \bw' = \cos^2(\chi/2) \bw + \sin^2(\chi/2) \bw_1 -
        |\bs| \cos(\chi/2)\sin(\chi/2) \bn_{\bw} = \bv, \\
      \end{gathered}
    \end{equation}
    Following \eqref{eq:conser_v}, and by change of variables, we arrive
    at
      \begin{equation}
        \label{eq:partA}
        \begin{aligned}
          &\int_{\mathbb{R}^3} \int_{\mathbb{R}^3} \int_{\bn \perp \bg}
          \int_0^{\pi} B(|\bg|,\chi)
          H^{\newindex{i}}(\bv') H^{\newindex{j}}(\bv'_1) H^{\indexk}(\bv) \exp
          \left( -\frac{|\bv|^2 + |\bv_1|^2}{2} \right) \,\mathrm{d}\chi
          \,\mathrm{d}\bn \,\mathrm{d}\bv_1 \,\mathrm{d}\bv  \\
          ={} & \int_{\mathbb{R}^3} \int_{\mathbb{R}^3} \int_{\tilde{\bn} \perp
            \bs} \int_0^{\pi} B(|\bs|,\chi)
          H^{\newindex{i}}(\bw) H^{\newindex{j}}(\bw_1) H^{\indexk}(\bw') \exp
          \left( -\frac{|\bw|^2 + |\bw_1|^2}{2} \right) \,\mathrm{d}\chi
          \,\mathrm{d}\tilde{\bn} \,\mathrm{d}\bw_1 \,\mathrm{d}\bw \\
          ={} &\int_{\mathbb{R}^3} \int_{\mathbb{R}^3} \int_{\bn \perp \bg}
          \int_0^{\pi} B(|\bg|,\chi)
          H^{\newindex{i}}(\bv) H^{\newindex{j}}(\bv_1) H^{\indexk}(\bv') \exp
          \left( -\frac{|\bv|^2 + |\bv_1|^2}{2} \right) \,\mathrm{d}\chi
          \,\mathrm{d}\bn \,\mathrm{d}\bv_1 \,\mathrm{d}\bv.
        \end{aligned}
      \end{equation}
    Thus, we can substitute the above equality into
    \eqref{eq:coeA_detail} to get
    \begin{equation}
      \label{eq:coeA_detail1}
      \begin{aligned}
        A_{k_1k_2k_3}^{\newindex{i}, \newindex{j}} =
        &\frac{1}{(2\pi)^3\factorialk} \int_{\mathbb{R}^3}
        \int_{\mathbb{R}^3} \int_{\bn \perp \bg} \int_0^{\pi} B(|\bg|,\chi)
        [H^{k_1 k_2 k_3}(\bv') - H^{k_1 k_2k_3}(\bv)] \\
        & \qquad H^{\newindex{i}}(\bv) H^{\newindex{j}}(\bv_1) \exp
        \left(
          -\frac{|\bv|^2 + |\bv_1|^2}{2} \right) \,\mathrm{d}\chi
        \,\mathrm{d}\bn \,\mathrm{d}\bv_1 \,\mathrm{d}\bv.
      \end{aligned}
    \end{equation}
    Further simplification of \eqref{eq:coeA_detail1} follows the
    method in \cite{Grad}, where the velocity of the mass center is
    defined as $\bh = (\bv+\bv_1)/2 = (\bv' + \bv'_1)/2$. Hence,
    \begin{gather}
      \label{eq:var_change}
      \bv = \bh + \frac{1}{2} \bg, \quad
      \bv_1 = \bh - \frac{1}{2} \bg, \quad
      \bv' = \bh + \frac{1}{2} \bg', \quad
      \bv_1' = \bh - \frac{1}{2} \bg', \\
      \label{eq:gh}
      |\bv|^2 + |\bv_1|^2 = \frac{1}{2} |\bg|^2 + 2|\bh|^2, \qquad
      \mathrm{d}\bv \, \mathrm{d}\bv_1 = \mathrm{d}\bg \, \mathrm{d}\bh.
    \end{gather}
    Combining Lemma \ref{thm:Hermit_cv}, Corollary \ref{thm:Hermit_split}
    and \eqref{eq:gh}, we can rewrite \eqref{eq:coeA_detail1} as an
    integral with respect to $\bg$ and $\bh$:
    \begin{equation}
      \label{eq:coeA_detail2}
      \begin{aligned}
        A_{k_1k_2k_3}^{\newindex{i}, \newindex{j}} = &
        \sum\limits_{i'_1+j'_1=i_1+j_1}\sum\limits_{i'_2+j'_2=i_2+j_2}\sum\limits_{i'_3+j'_3=i_3+j_3}
        \sum\limits_{l'_1+k'_1=k_1}\sum\limits_{l'_2+k'_2=k_2}\sum\limits_{l'_3+k'_3=k_3}
        \\ & \qquad \frac{
          2^{-k/2}}{(2\pi)^3}\frac{1}{k_1'k_2'k_3'l'_1!l'_2!l'_3!}
        a_{i'_1j'_1}^{i_1j_1}a_{i'_2j'_2}^{i_2j_2}
        a_{i'_3j'_3}^{i_3j_3}\gamma_{j'_1j'_2j'_3}^{l'_1l'_2l'_3}\eta_{i'_1i'_2i'_3}^{k'_1k'_2k'_3},
      \end{aligned}
    \end{equation}
    where the coefficients $\gamma_{j'_1j'_2j'_3}^{l'_1l'_2l'_3}$
    defined in \eqref{eq:coe_gamma} are integrals with respect to
    $\bg$, and $\eta_{i'_1 i'_2 i'_3}^{k'_1 k'_2 k'_3}$ are integrals
    with respect to $\bh$ defined by
    \begin{equation}
      \label{eq:eta}
      \eta_{i'_1 i'_2 i'_3}^{k'_1 k'_2 k'_3} = \int_{\mathbb{R}^3}
      H^{i'_1 i'_2 i'_3} (\sqrt{2} \bh) H^{k'_1 k'_2 k'_3} (\sqrt{2} \bh)
      \exp (-|\bh|^2) \,\mathrm{d}\bh =
      \pi^{3/2} k'_1! k'_2! k'_3! \delta_{i'_1 k'_1} \delta_{i'_2 k'_2}
      \delta_{i'_3 k'_3}.
    \end{equation}
    Thus the theorem is proven by substituting \eqref{eq:eta} into
    \eqref{eq:coeA_detail2}.
  \end{proof}
}

\section{Proof of Theorem \ref{thm:gamma}}
We will first prove Theorem \ref{thm:gamma} based on
several lemmas, and then prove these lemmas. 
\subsection {Proof of Theorem \ref{thm:gamma}} 
In order to prove Theorem \ref{thm:gamma}, we will introduce the
definition of Ikenberry polynomials \cite{Ikenberry1955} and several
lemmas.
\begin{definition}[Ikenberry polynomials] \label{def:hhp}
  Let $\bv = (v_1, v_2, v_3)^T \in \mathbb{R}^3$. For
  $\forall n \in \bbN$, and $i_1, \cdots, i_n \in \{1,2,3\}$, define
  $Y_{i_1 \cdots i_n}(\bv)$ as the Ikenberry polynomials 
\begin{gather*}
  Y(\bv) = 1, \qquad Y_{i_1}(\bv) =  v_{i_1}, \\
  Y_{i_1 \cdots i_n}(\bv) = v_{i_1}\cdots v_{i_n} +
  |\bv|^2S_{n-2}^{i_1\cdots i_n}(\bv) + |\bv|^4S_{n-4}^{i_1\cdots i_n}
  + \cdots + |\bv|^{2\lfloor n/2 \rfloor}S_{n-2\lfloor n/2
    \rfloor}^{i_1\cdots i_n}(\bv),
\end{gather*}
where $S_{j}^{i_1\cdots i_n}$ is a homogeneous harmonic polynomial of
degree $j$ defined in \cite{Ikenberry1955}, which can be determined by
\begin{equation*}
  \Delta_{\bv}Y_{i_1\cdots i_n} = \Delta_{\bv}^2Y_{i_1\cdots i_n} =
  \Delta_{\bv}^{\lfloor n/2 \rfloor}Y_{i_1\cdots i_n} = 0. 
\end{equation*}
For $k_1, k_2, k_3 \in \mathbb{N}$, define $Y^{k_1 k_2 k_3}(\bv)$ as
the polynomial $Y_{i_1 \cdots i_n}(\bv)$ with
\begin{gather*}
n = k_1 + k_2 + k_3, \quad i_1 = \cdots = i_{k_1} = 1, \\
i_{k_1+1} = \cdots = i_{k_1+k_2} = 2, \quad
i_{k_1+k_2+1} = \cdots = i_n = 3.
\end{gather*}
\end{definition}
\begin{lemma}
\label{thm:S_harmonic}
The integral
\begin{displaymath}
\int_{\mathbb{S}^2} Y^{k_1 k_2 k_3}(\bn) Y^{l_1 l_2 l_3}(\bn) \,\mathrm{d}\bn
\end{displaymath}
is the coefficient of $v_1^{k_1} v_2^{k_2} v_3^{k_3} w_1^{l_1}
w_2^{l_2} w_3^{l_3}$ in the polynomial
\begin{displaymath}
\frac{4\pi}{2k+1} \frac{k_1! k_2! k_3! l_1! l_2! l_3!}{[(2k-1)!!]^2}
(|\bv| |\bw|)^k P_k \left(
  \frac{\bv}{|\bv|} \cdot \frac{\bw}{|\bw|}
\right), \qquad k = k_1 + k_2 + k_3.
\end{displaymath}
\end{lemma}
\begin{lemma}
  \label{thm:Hermite_Harmonic}
  The Hermite polynomial $H^{k_1 k_2 k_3}(\bv)$ can be represented as
  \begin{displaymath}
    \begin{split}
      H^{k_1 k_2 k_3}(\bv) &=
      \sum_{m_1=0}^{\lfloor k_1/2 \rfloor}
      \sum_{m_2=0}^{\lfloor k_2/2 \rfloor}
      \sum_{m_3=0}^{\lfloor k_3/2 \rfloor}
      \frac{(-1)^m m! (2k-4m+1)!!} {(2(k-m)+1)!!}
      \left( \prod_{i=1}^3 \frac{k_i!}{m_i! (k_i-2m_i)!} \right) \\
      &\qquad L_m^{(k-2m+1/2)} \left( \frac{|\bv|^2}{2} \right) 
      Y^{k_1-2m_1,k_2-2m_2,k_3-2m_3}(\bv),
    \end{split}
  \end{displaymath}
  where $k = k_1 + k_2 + k_3$ and $m = m_1 + m_2 + m_3$.
\end{lemma}

\begin{lemma}
  \label{thm:int_n}
  Given a vector $\bg$ and $\chi \in [0, \pi]$, let
  $\bg'(\bn) = \bg \cos \chi - |\bg| \bn \sin \chi$, where $\bn$ is a unit
  vector. It holds that
  \begin{displaymath}
    \int_{\bn \perp \bg} Y^{k_1 k_2 k_3}(\bg' / |\bg|) \,\mathrm{d}\bn
    = 2\pi Y^{k_1 k_2 k_3}(\bg / |\bg|) P_k(\cos \chi),
  \end{displaymath}
  where $k = k_1 + k_2 + k_3$ and $P_k$ is Legendre polynomial.
\end{lemma}

In above lemmas, Lemma \ref{thm:S_harmonic} and Lemma
\ref{thm:Hermite_Harmonic} will be proved in Appendix
\ref{sec:S_harmonic} and \ref{sec:Hermite_Harmonic}
respectively. Lemma \ref{thm:int_n} is proved in \cite{Ikenberry1956}.
By Lemma \ref{thm:Hermite_Harmonic} and Lemma
\ref{thm:int_n}, we can derive the corollary below
\begin{corollary}
  \label{thm:int_H}
  Given a vector $\bg$ and $\chi \in [0, \pi]$, define $\bg'(\bn)$ the
  same as in Theorem \ref{thm:int_n}. We have
  \begin{displaymath}
    \begin{split}
      & \int_{\bn \perp \bg} H^{k_1 k_2 k_3}(\bg') \,\mathrm{d}\bn =
      2\pi \sum_{m_1=0}^{\lfloor k_1/2 \rfloor}
      \sum_{m_2=0}^{\lfloor k_2/2 \rfloor}
      \sum_{m_3=0}^{\lfloor k_3/2 \rfloor}
      \frac{(-1)^m m! (2k-4m+1)!!} {(2(k-m)+1)!!} \times \\
      & \quad \left( \prod_{i=1}^3 \frac{k_i!}{m_i! (k_i-2m_i)!} \right)
      L_n^{(k-2m+1/2)} \left( \frac{|\bg|^2}{2} \right) 
      Y^{k_1-2m_1,k_2-2m_2,k_3-2m_3}(\bg) P_{k-2m}(\cos \chi),
    \end{split}
  \end{displaymath}
  where $k = k_1 + k_2 + k_3$, $m = m_1 + m_2 + m_3$.
\end{corollary}

{\renewcommand\proofname{Proof of Theorem \ref{thm:gamma}}
\begin{proof}

  By Lemma \ref{thm:Hermite_Harmonic}, the corollary \ref{thm:int_H}
  and the homogeneity of the Ikenberry polynomials 
  $ \gamma_{\indexk}^{\indexl}$ defined in \eqref{eq:gamma} can
    be simplified as
{\small 
  \begin{equation}
    \label{eq:gamma_detail2}
    \begin{aligned}
      \gamma_{\indexk}^{\indexl} & = 2\pi \sum_{m_1=0}^{\lfloor k_1/2
        \rfloor} \sum_{m_2=0}^{\lfloor k_2/2 \rfloor}
      \sum_{m_3=0}^{\lfloor k_3/2 \rfloor}\sum_{n_1=0}^{\lfloor l_1/2
        \rfloor} \sum_{n_2=0}^{\lfloor l_2/2 \rfloor}
      \sum_{n_3=0}^{\lfloor l_3/2 \rfloor}
      \frac{(2(k-m)+1)!!C_{\indexm}^{\indexk}}{4\pi\prod_{i=1}^3
        (k_i-2m_i)!}\frac{(2(l-n)+1)!!C_{\indexn}^{\indexl}}{4\pi\prod_{i=1}^3
        (l_i-2n_i)!} \times \\
      & \qquad 2\int_{0}^{+\infty} \int_0^{\pi}\int_{\bbS^2}
      Y^{k_1-2m_1,k_2-2m_2,k_3-2m_3}(\bn)
      Y^{l_1-2n_1,l_2-2n_2,l_3-2n_3}(\bn)\left( \frac{g}{\sqrt{2}} \right)^{k+l+2-2(m+n)} \times \\
      & \qquad \qquad L_m^{(k-2m+1/2)}\left( \frac{g^2}{4} \right)L_n^{(l-2n+1/2)}
      \left( \frac{g^2}{4} \right) B(g, \chi)\Big[P_{k-2m}(\cos \chi)- 1\Big] \exp\left(
        -\frac{g^2}{4} \right) \dd \bn \,\mathrm{d}\chi\,\mathrm{d}g,
    \end{aligned}
  \end{equation}
}
where $C_{\indexm}^{\indexl}$ is defined in \eqref{eq:coeC}.

Substituting Lemma \ref{thm:S_harmonic} into \eqref{eq:gamma_detail2},
we complete this proof.
\end{proof}
}

\subsection{Proof of Lemma \ref{thm:S_harmonic}}
\label{sec:S_harmonic}
In order to prove Lemma \ref{thm:S_harmonic}, we first introduce the
following definitions and lemmas.
\begin{definition}[Associated Legendre functions]
  For $m = -l,\cdots,l$, the associated Legendre functions are
  defined as
  \begin{displaymath}
    P_{l}^m(x) = \frac{(-1)^m}{2^{l} l!} (1-x^2)^{m/2}
    \frac{\mathrm{d}^{l+m}}{\mathrm{d}x^{l+m}} (x^2-1)^{l}.
  \end{displaymath}
\end{definition}
\begin{definition}[Spherical harmonics]
  For $l \in \mathbb{N}$ and $m = -l, \cdots, l$, the spherical
  harmonic $Y_{l}^m(\theta, \varphi)$ is defined as
  \begin{displaymath}
    Y_{l}^m(\bn) = Y_{l}^m(\theta, \varphi) =
    \sqrt{\frac{2l+1}{4\pi} \frac{(l - m)!}{(l + m)!}}
    P_{l}^m(\cos \theta) \exp(\mathrm{i} m \varphi),
    \qquad \bn \in \mathbb{S}^2,
  \end{displaymath}
  where $(\theta,\varphi)$ is the spherical coordinates of $\bn$.
\end{definition}

\begin{lemma}[Addition theorem]
  \label{thm:addition_theorem}
  For any $l \in \mathbb{N}$, it holds that
  \begin{displaymath}
    P_{l}(\bn_1 \cdot \bn_2) = \frac{4\pi}{2l+1}
    \sum_{m=-l}^{l} Y_{l}^m(\bn_1)
    \overline{Y_{l}^m(\bn_2)},
  \end{displaymath}
  where $P_{l}$ is  Legendre polynomial.
\end{lemma}
\begin{lemma} 
  \label{thm:P=Y} 
  For any $l \in \mathbb{N}$, it holds that
  \begin{displaymath}
    (|\bv| |\bw|)^{l} P_{l} \left(
      \frac{\bv}{|\bv|} \cdot \frac{\bw}{|\bw|}
    \right) = \frac{(2l)!}{2^{l} l! l!}
    \sum_{i_1=1}^3 \cdots \sum_{i_{l}=1}^3
    w_{i_1}\cdots w_{i_{l}}
    Y_{i_1\cdots i_{l}}(\bv).
  \end{displaymath}
\end{lemma}

In the above lemmas, Lemma \ref{thm:addition_theorem} and Lemma
\ref{thm:P=Y} are well-known and their proofs can be found in
\cite{1999AmJPh} and \cite{Ikenberry1961} respectively. Based on these
two lemmas, the following corollary holds.
\begin{corollary}
  \label{thm:coe_Y}
  The harmonic polynomial $Y^{k_1 k_2 k_3}(\bv)$ is the coefficient of
  the monomial $w_1^{k_1} w_2^{k_2} w_3^{k_3}$ in the following
  polynomial of $\bw$:
  \begin{displaymath}
    \frac{k_1! k_2! k_3!}{(2k-1)!!} (|\bv| |\bw|)^k P_k \left(
      \frac{\bv}{|\bv|} \cdot \frac{\bw}{|\bw|}
    \right), \qquad k = k_1 + k_2 + k_3.
  \end{displaymath}
\end{corollary}
{\renewcommand\proofname{Proof of Corollary \ref{thm:coe_Y}}
  \begin{proof}
    Since 
    \begin{displaymath}
      \sum_{i_1=1}^3 \cdots \sum_{i_{k}=1}^3
      w_{i_1}\cdots w_{i_{k}}
      Y_{i_1\cdots i_{k}}(\bv) =  \frac{k!}{k_1!k_2!k_3!}
      \sum_{k_1+k_2+k_3=k} w_1^{k_1}w_2^{k_2}w_3^{k_3}
      Y^{k_1k_2k_3}(\bv),
    \end{displaymath}
    and matching the term of $w_1^{k_1}w_2^{k_2}w_3^{k_3}$ in Lemma
    \ref{thm:P=Y}, we complete this proof.
  \end{proof}
}

{\renewcommand\proofname{Proof of Lemma \ref{thm:S_harmonic}}
  \begin{proof}
    From Corollary \ref{thm:coe_Y}, we can derive that
    $\int_{\mathbb{S}^2} Y^{k_1 k_2 k_3}(\bn) Y^{l_1 l_2 l_3}(\bn)
    \,\mathrm{d}\bn$
    is the coefficient of
    $v_1^{k_1}v_2^{k_2}v_3^{k_3}w_1^{l_1}w_2^{l_2} w_3^{l_3}$ in the polynomial
    \begin{displaymath}
      \int_{\mathbb{S}^2} 
      \left[\beta^{k_1k_2k_3} (|\bn| |\bv|)^kP_k\left(\bn \cdot
          \frac{\bv}{|\bv|}\right) \right]\left[\beta^{l_1l_2l_3} (|\bn|
        |\bw|)^{l}P_{l}\left(\bn \cdot
          \frac{\bw}{|\bw|}\right) \right]
      \dd \bn,
    \end{displaymath}
    where $k = k_1 + k_2 + k_3, l = l_1 + l_2 + l_3$ and
    $\beta^{\indexk} = \frac{\factorialk}{(2k-1)!!}$.  Following Theorem
    \ref{thm:addition_theorem}, it holds 
    \begin{align*}
      &   \int_{\mathbb{S}^2} 
        \left[ (|\bn| |\bv|)^kP_k\left(\bn \cdot
        \frac{\bv}{|\bv|}\right) \right]\left[ (|\bn|
        |\bw|)^{l}P_{l}\left(\bn \cdot
        \frac{\bw}{|\bw|}\right) \right]
        \dd \bn \\
      & =(|\bv|^k |\bw|)^l
        \frac{(4\pi)^2}{(2k+1)(2l+1)}
        \sum_{m=-k}^{k} \sum_{n=-l}^{l}  Y_{k}^m(\bv)
        \overline{Y_{l}^n(\bw)}\delta_{lk}\delta_{mn} \\
      & = \frac{4\pi\delta_{kl}}{2k+1}(|\bv| |\bw|)^k
        P_k\left(\frac{\bv}{|\bv|} \cdot
        \frac{\bw}{|\bw|}\right).
    \end{align*}
    Thus if $k = l$, this corollary is proved. If $k \neq l$, we can
    deduce that
    $\int_{\mathbb{S}^2} Y^{k_1 k_2 k_3}(\bn) Y^{l_1 l_2 l_3}(\bn)
    \,\mathrm{d}\bn = 0$.
    In this case, the coefficient of
    $v_1^{k_1}v_2^{k_2}v_3^{k_3}w_1^{l_1}w_2^{l_2} w_3^{l_3}$ in the
    polynomial
    $(|\bv| |\bw|)^k P_k\left(\frac{\bv}{|\bv|} \cdot
      \frac{\bw}{|\bw|}\right)$
    is also zero, and this completes the proof.
  \end{proof}
}

\subsection{Proof of Lemma \ref{thm:Hermite_Harmonic}}
\label{sec:Hermite_Harmonic}
We will prove Lemma \ref{thm:Hermite_Harmonic} in this section.  
{\renewcommand\proofname{Proof of Lemma \ref{thm:Hermite_Harmonic}}
  \begin{proof}
    Define the homogeneous spherical harmonic
    $Z_{i_1i_2\cdots i_k}^{(k,m)}$ of degree $k-2m$ as
    \begin{equation}
      \label{eq:hsh_Z}
      Z_{i_1i_2\cdots i_k}^{(k,m)}
        = \frac{1}{k!} \sum_{\sigma \in \mathcal{S}_k}
          Y_{i_{\sigma(1)} i_{\sigma(2)} \cdots i_{\sigma(r)}}
          \delta_{i_{\sigma(r+1)} i_{\sigma(r+2)}} \cdots
	  \delta_{i_{\sigma(k-1)} i_{\sigma(k)}},  
    \end{equation}
    where $r = k - 2m$ and the sum is taken over all permutations of
    the set $\{1,2,\cdots,k\}$, i.e.
    \begin{displaymath}
    \mathcal{S}_k = \{ \sigma \mid \sigma: \{1,2,\cdots,k\}
      \rightarrow \{1,2,\cdots,k\} \text{ is a bijection} \}.
    \end{displaymath}
    It has been proven in \cite[eqs. (3)(8)(9)(31)]{Ikenberry1962} that 
    \footnote{In \cite{Ikenberry1962}, the definition of the Laguerre
    polynomial differs from Definition \ref{def:Laguerre} by a
    constant, which makes the coefficient in our paper slightly
    different from the one in \cite{Ikenberry1962}.}
    \begin{equation}
      \label{eq:Hermit_split}
      H^{k_1k_2k_3}(\bv) =
      \sum_{m=0}^{\lfloor k/2 \rfloor}
      \frac{(-1)^mk!(2k - 4m +1)!!}{(k-2m)!(2k-2m+1)!!}
      L_m^{(k-2m+1/2)} \left( \frac{|\bv|^2}{2} \right) 
      Z^{(k,m)}_{i_1i_2\cdots i_k}(\bv),
    \end{equation}
    where the indices $i_1, \cdots, i_k$ satisfy:
    \begin{displaymath}
    i_1 = \cdots = i_{k_1} = 1, \qquad
    i_{k_1+1} = \cdots = i_{k_1+k_2} = 2, \qquad
    i_{k_1+k_2+1} = \cdots = i_k = 3.
    \end{displaymath}

    To prove Lemma \ref{thm:Hermite_Harmonic}, we just need to provide
    a more explicit expression for \eqref{eq:hsh_Z}. In order that the
    summand in \eqref{eq:hsh_Z} is nonzero, the two indices of every
    Kronecker symbol must be the same. When all the Kronecker symbols
    take $2m_1$ ones, $2m_2$ twos and $2m_3$ threes as their indices,
    the summand will actually be $Y^{k_1-2m_1, k_2-2m_2,
    k_3-2m_3}(\bv)$ according to Definition \ref{def:hhp}. Apparently
    $(m_1, m_2, m_3)$ must be indices from the following set:
    \begin{displaymath}
    \mathcal{M}_{k_1k_2k_3}^m = \{(m_1, m_2, m_3) \mid
      m_1 + m_2 + m_3 = m, \, 2m_1 \leqslant k_1, \,
      2m_2 \leqslant k_2, \, 2m_3 \leqslant k_3 \}.
    \end{displaymath}

    Next, we are going to count how many times $Y^{k_1-2m_1, k_2-2m_2,
    k_3-2m_3}(\bv)$ appears in the sum in \eqref{eq:hsh_Z}. This can
    be observed by noting that
    \begin{enumerate}
    \item The $m$ Kronecker symbols choosing from $m_1$ pairs of ones,
      $m_2$ pairs of twos and $m_3$ pairs of threes gives a factor
      $m! / (m_1! m_2! m_3!)$;
    \item The $k-2m$ indices of $Y$ choosing from $k_1-2m_1$ ones,
      $k_2-2m_2$ twos and $k_3-2m_3$ threes gives a factor $(k-2m)! /
      \big( (k_1-2m_1)! (k_2-2m_2)! (k_3-2m_3)! \big)$.
    \item Permutations of $k_1$ ones, $k_2$ twos and $k_3$ threes give
      respectively factors $k_1!$, $k_2!$ and $k_3!$.
    \end{enumerate}
    Summarizing all these results, we get
    \begin{equation}
      \label{eq:Z}
      \begin{split}
        & Z^{(k,m)}_{i_1i_2\cdots i_k} = \frac{1}{k!}
	  \sum_{(m_1,m_2,m_3)\in \mathcal{M}_{k_1k_2k_3}^m}
	    \frac{(k-2m)!m!\prod\limits_{i=1}^3 k_i!}%
	      {\prod\limits_{i=1}^3 \Big( (k_i-2m_i)!m_i!\Big)}
	    Y^{k_1-2m_1, k_2-2m_2, k_3-2m_3}(\bv).
      \end{split}
    \end{equation}
    By \eqref{eq:Hermit_split} and \eqref{eq:Z}, the proof is
    completed.
  \end{proof}
}

\bibliographystyle{plain}
\bibliography{article}
\end{document}